\documentclass{cedram-alco}
% ...
\usepackage{amsmath}
\usepackage{amsthm}
\usepackage{amssymb}
\usepackage{amscd}
\usepackage{color,graphicx}
\usepackage[all]{xy}
\usepackage{caption,subcaption}
\newcommand\vectz{\boldsymbol{z}}
\newcommand\vecta{\boldsymbol{a}}
\newcommand\Nat{\mathbb{N}}
\newcommand\Int{\mathbb{Z}}
\newcommand\Rat{\mathbb{Q}}
\newcommand\Comp{\mathbb{C}}
\renewcommand\tilde{\widetilde}

\newcommand\Hom{\operatorname{Hom}}

\newcommand\AP{\mathcal{A}}
\newcommand\BB{\mathcal{B}}
\newcommand\EE{\mathcal{E}}
\newcommand\GG{\mathcal{G}}
\newcommand\LL{\mathcal{L}}
\newcommand\OO{\mathcal{O}}
\newcommand\PP{\mathcal{P}}
\newcommand\RR{\mathcal{R}}
\newcommand\TT{\mathcal{T}}
\newcommand\XX{\mathcal{X}}
\newcommand\point{\operatorname{pt}}

\title
{Skew hook formula for $d$-complete posets via equivariant $K$-theory}

\author[\initial{H.} Naruse]{\firstname{Hiroshi} \lastname{Naruse}}

\address{%
Hiroshi Naruse\\
Graduate School of Education, University of Yamanashi\\ 
4-4-37, Takeda, Kofu, Yamanashi 400-8510, Japan}

\email{hnaruse@yamanashi.ac.jp}

\author[\initial{S.} Okada]{\firstname{Soichi} \lastname{Okada}}

\address{%
Soichi Okada\\
Graduate School of Mathematics, Nagoya University\\ 
Furo-cho, Chikusa-ku, Nagoya 464-8602, Japan}

\email{okada@math.nagoya-u.ac.jp}

\keywords{%
$d$-complete posets, hook formulas, $P$-partitions, Schubert calculus, equivariant $K$-theory
}

\subjclass{%
05A15 (primary), 06A07, 14N15, 19L47 (secondary)
}

\begin{document}
\begin{abstract}
Peterson and Proctor obtained a formula which expresses the multivariate generating function 
for $P$-partitions on a $d$-complete poset $P$ as a product in terms of hooks in $P$.
In this paper, we give a skew generalization of Peterson--Proctor's hook formula, 
\ie a formula for the generating function of $(P \setminus F)$-partitions for 
a $d$-complete poset $P$ and an order filter $F$ of $P$, by using the notion of excited diagrams.
Our proof uses the Billey-type formula and the Chevalley-type formula 
in the equivariant $K$-theory of Kac--Moody partial flag varieties. 
This generalization provides an alternate proof of Peterson--Proctor's hook formula.
As the equivariant cohomology version, we derive a skew generalization 
of a combinatorial reformulation of Nakada's colored hook formula for roots.
\end{abstract}

\maketitle

\section{%
Introduction
}

One of the most elegant formulas in combinatorics is the Frame--Robinson--Thrall hook formula 
\cite[Theorem~1]{FRT} for the number of standard tableaux.
Given a partition $\lambda$ of $n$, 
a standard tableaux of shape $\lambda$ is a filling of the cells of 
the Young diagram $D(\lambda)$ of $\lambda$ with numbers $1, 2, \dots, n$ 
such that each number appears once and the entries of each row and each column are increasing.
The Frame--Robinson--Thrall hook formula asserts that the number $f^\lambda$ of 
standard tableaux of shape $\lambda$ is given by
\begin{equation}
\label{eq:FRT-hook}
f^\lambda
 =
\frac{ n! }
     { \prod_{v \in D(\lambda)} h_{D(\lambda)}(v) },
\end{equation}
where $h_{D(\lambda)}(v)$ denotes the hook length of the cell $v$ in $D(\lambda)$.
Similar formulas hold for the number of shifted standard tableaux 
(\cite[5.1.4, Exercise~21]{Knu}, see also \cite[\S~40]{Schur} and \cite[Theorem~1]{T}) 
and the number of increasing labeling of rooted trees (\cite[5.1.4, Exercise~20]{Knu}).
These tableaux and labelings can be regarded as linear extensions of certain posets.

Stanley \cite{Stan1} introduced the notion of $P$-partitions for a poset $P$, 
and found a relationship between the univariate generating function 
and the number of linear extensions of $P$.
Given a poset $P$, a $P$-partition is an order-reversing map $\sigma$ from $P$ to $\Nat$, 
the set of nonnegative integers.
We denote by $\AP(P)$ the set of all $P$-partitions.
For a $P$-partition $\sigma$, we write $|\sigma| = \sum_{v \in P} \sigma(x)$.
Then Stanley \cite[Corollaries~5.3 and 5.4]{Stan1} proved that, for a poset $P$ with $n$ elements, 
there exists a polynomial $W_P(q)$ satisfying
\begin{equation}
\label{eq:P-partition}
\sum_{\sigma \in \AP(P)} q^{|\sigma|}
 =
\frac{ W_P(q) }
     { \prod_{i=1}^n (1 - q^i) },
\end{equation}
and that $W_P(1)$ is equal to the number of linear extensions of $P$.
Also in \cite[Proposition~18.3]{Stan2} he proved that, 
if $P$ is the Young diagram $D(\lambda)$ of a partition $\lambda$, viewed as a poset, 
the generating function of $D(\lambda)$-partitions (also called reverse plane partitions of shape $\lambda$) 
is given by
\begin{equation}
\label{eq:Stanley-hook}
\sum_{\sigma \in \AP(D(\lambda))} q^{|\sigma|}
 =
\frac{ 1 }
     { \prod_{v \in D(\lambda)} (1 - q^{h_{D(\lambda)}(v)}) }.
\end{equation}
Combining \eqref{eq:P-partition} and \eqref{eq:Stanley-hook} and taking the limit $q \to 1$, 
we obtain the Frame--Robinson--Thrall hook formula \eqref{eq:FRT-hook}.
Gansner \cite[Theorem~5.1]{G} gave a multivariate generalization of \eqref{eq:Stanley-hook}.

Proctor \cite{Proc1}, \cite{Proc2} introduced a wide class of posets, 
called $d$-complete posets, enjoying ``hook-length property'', 
as a generalization of Young diagrams, shifted Young diagrams and rooted trees.
$d$-Complete posets are defined by certain local structural conditions 
(see Section~2 for a precise definition).
Peterson and Proctor obtained the following theorem, which is a far-reaching generalization 
of the hook formulas \eqref{eq:FRT-hook} and \eqref{eq:Stanley-hook}.

\begin{theo}
\label{thm:PP}
\textup{(Peterson--Proctor, see \cite{Proc3})}
Let $P$ be a $d$-complete poset.
The multivariate generating function of $P$-partitions is given by
\begin{equation}
\label{eq:PP-hook}
\sum_{\sigma \in \AP(P)} \vectz^{\sigma}
 =
\frac{ 1 }
     { \prod_{v \in P} (1 - \vectz[H_P(v)]) }.
\end{equation}
\textup{(}Refer to Section~2 for undefined notations.\textup{)}
\end{theo}

However the original proof of this theorem is not yet published, 
though an outline of their proof is given in \cite{Proc3}.
Different proofs are sketched by Ishikawa--Tagawa \cite{IT1}, \cite{IT2} and Nakada \cite{N2}, \cite{N3}.
Our skew generalization (Theorem~\ref{thm:main} below) provides an alternate proof of Theorem~\ref{thm:PP}.
In the univariate case, a full proof is given by Kim--Yoo \cite{KY}.

Another direction of generalizing the Frame--Robinson--Thrall hook formula \eqref{eq:FRT-hook} 
is to consider skew shapes.
For partitions $\lambda \supset \mu$, 
a standard tableau of skew shape $\lambda/\mu$ is a filling of the cells of the skew Young diagram 
$D(\lambda/\mu) = D(\lambda) \setminus D(\mu)$ satisfying the same conditions as standard tableaux 
of straight shape.
However one cannot expect a nice product formula for the number $f^{\lambda/\mu}$ of 
standard tableaux of skew shape $\lambda/\mu$ in general.
Naruse \cite{Naruse} presented and sketched a proof of a subtraction-free formula for $f^{\lambda/\mu}$:
\begin{equation}
\label{eq:Naruse-hook}
f^{\lambda/\mu}
 =
n!
\sum_{D \in \EE_{D(\lambda)}(D(\mu))} 
 \frac{ 1 }
      { \prod_{v \in D(\lambda) \setminus D} h_\lambda(v) },
\end{equation}
where $n = |\lambda/\mu|$, 
and $D$ runs over all excited diagrams of $D(\mu)$ in $D(\lambda)$.
Morales--Pak--Panova \cite{MPP} gave a $q$-analogue of Naruse's skew hook formula 
for the univariate generating functions for $P$-partitions on $P = D(\lambda/\mu)$.

The main result of this paper is the following skew generalization of Peterson--Proctor's 
hook formula (Theorem~\ref{thm:PP}).
Recall that a subset $F$ of a poset $P$ is called an order filter of $P$ if 
$x<y$ in $P$ and $x \in F$ imply $y \in F$.
In particular the empty set $F = \emptyset$ is an order filter of $P$.

\begin{theo}
\label{thm:main}
Let $P$ be a $d$-complete poset and $F$ an order filter of $P$.
Then the multivariate generating function of $(P \setminus F)$-partitions, 
where $P \setminus F$ is viewed as an induced subposet of $P$, is given by
\begin{equation}
\label{eq:main}
\sum_{\sigma \in \AP(P \setminus F)} \vectz^\sigma
 =
\sum_{D \in \EE_P(F)}
 \frac{ \prod_{v \in B(D)} \vectz[H_P(v)] }
      { \prod_{v \in P \setminus D} ( 1 - \vectz[H_P(v)] ) },
\end{equation}
where $D$ runs over all excited diagrams of $F$ in $P$.
\textup{(}See Sections~2 and 3 for undefined notations.\textup{)}
\end{theo}

Taking an appropriate limit, we see that 
the number of linear extensions of $P \setminus F$ is given by
\begin{equation}
\label{eq:lin_ext}
n! 
\sum_{D \in \EE_P(F)}
 \frac{ 1 }
      { \prod_{v \in P \setminus D} h_P(v) },
\end{equation}
where $n = \# (P \setminus F)$ and $h_P(v)$ is the hook length of $v$ in $P$.
(See Corollary~\ref{cor:main3}(b).)

If $F = \emptyset$, then our main theorem (Theorem~\ref{thm:main}) gives 
Peterson--Proctor's hook formula (Theorem~\ref{thm:PP}).
If $P = D(\lambda)$ and $F = D(\mu)$ are the Young diagrams of partitions $\lambda \supset \mu$, 
then \eqref{eq:main} reduces to Morales--Pak--Panova's $q$-hook formula \cite[Corollary~6.17]{MPP} 
after specializing $z_i = q$ for all $i \in I$, 
and \eqref{eq:lin_ext} is nothing but Naruse's skew hook formula \eqref{eq:Naruse-hook}.

Our proof of Theorem~\ref{thm:main} uses the equivariant $K$-theory of Kac--Moody partial flag varieties.
In fact, the notion of excited diagrams and excited peaks for ordinary and shifted Young diagrams 
has its origin in the equivariant Schubert calculus 
(see \cite{IN1}, \cite{IN2}; see also \cite{GK}, \cite{KN}, \cite{Kr1}, \cite{Kr2}).
Theorem~\ref{thm:main} is obtained by proving that 
the both sides of \eqref{eq:main} equal to the same ratio of the localized equivariant Schubert classes. 
(See Theorems~\ref{thm:xi1} and \ref{thm:xi2}.)
A key role is played by the Billey-type formula due to Lam--Schilling--Shimozono \cite{LSS} 
and the Chevalley-type formula due to Lenart--Shimozono \cite{LS} (see also \cite{LP}).

This paper is organized as follows.
In Section~2, we review a definition and basic properties of $d$-complete posets.
In Section~3, we introduce the notion of excited diagrams for $d$-complete posets, 
which is the key ingredient of the formulation of our main theorem, and study their properties.
In Section~4, we recall some properties of the equivariant $K$-theory 
and translate the Billey-type formula and the Chevalley-type formula in terms of combinatorics 
of $d$-complete posets.
We will give a proof of our main theorem (Theorem~\ref{thm:main}) and derive some corollaries 
in Section~5.

\section{%
$d$-Complete posets
}

In this section we review a definition and some properties of $d$-complete posets 
and explain their connections to Weyl groups.
See \cite{Proc1}, \cite{Proc2}, \cite{Proc3} and \cite{Stem2} for details.
We use the terminology for posets (partially ordered sets) 
from \cite[Chapter~3]{Stan3}.

\subsection{%
Combinatorics of $d$-complete posets
}

\begin{figure}[!htb]
\centering
\setlength{\unitlength}{1.5pt}
\begin{picture}(50,90)(-5,0)
\put(20,5){\circle*{3}}
\put(20,15){\circle*{3}}
\put(20,35){\circle*{3}}
\put(10,45){\circle*{3}}
\put(30,45){\circle*{3}}
\put(20,55){\circle*{3}}
\put(20,75){\circle*{3}}
\put(20,85){\circle*{3}}
\put(20,5){\line(0,1){10}}
\put(20,15){\line(0,1){5}}
\multiput(20,20)(0,2){5}{\line(0,1){1}}
\put(20,30){\line(0,1){5}}
\put(20,35){\line(-1,1){10}}
\put(20,35){\line(1,1){10}}
\put(10,45){\line(1,1){10}}
\put(30,45){\line(-1,1){10}}
\put(20,55){\line(0,1){5}}
\multiput(20,60)(0,2){5}{\line(0,1){1}}
\put(20,70){\line(0,1){5}}
\put(20,75){\line(0,1){10}}
\put(20,80){\makebox(15,10){$u_1$}}
\put(20,70){\makebox(15,10){$u_2$}}
\put(20,50){\makebox(25,10){$u_{k-2}$}}
\put(-5,40){\makebox(15,10){$x$}}
\put(30,40){\makebox(15,10){$y$}}
\put(20,30){\makebox(25,10){$v_{k-2}$}}
\put(20,10){\makebox(15,10){$v_2$}}
\put(20,0){\makebox(15,10){$v_1$}}
\end{picture}
\caption{Double-tailed diamond $d_k(1)$}
\label{fig:diamond}
\end{figure}
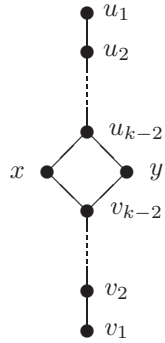

For an integer $k \ge 3$, we denote by $d_k(1)$ the poset consisting of $2k-2$ elements 
$u_1, \cdots, u_{k-2}, x, y, v_{k-2}, \cdots, v_1$ with covering relations
\begin{gather*}
u_1 \gtrdot u_2 \gtrdot \cdots \gtrdot u_{k-2}, \\
u_{k-2} \gtrdot x \gtrdot v_{k-2}, 
\quad u_{k-2} \gtrdot y \gtrdot v_{k-2}, \\
v_{k-2} \gtrdot \cdots \gtrdot v_2 \gtrdot v_1.
\end{gather*}
Note that $x$ and $y$ are incomparable.
The poset $d_k(1)$ is called the \emph{double-tailed diamond}.
The Hasse diagram of $d_k(1)$ is shown in Figure~\ref{fig:diamond}.

Let $P$ be a poset.
An interval $[v,u] = \{ x \in P : v \le x \le u \}$ is called 
a \emph{$d_k$-interval} if it is isomorphic to $d_k(1)$.
Then $v$ and $u$ are called the \emph{bottom} and \emph{top} of $[v,u]$ respectively,
and the two incomparable elements of $[v,u]$ are called the \emph{sides}.
A subset $I$ of $P$ is called \emph{convex} if $x < y < z$ in $P$ and $x$, $z \in I$ imply $y \in I$.
A convex subset $I$ is called a \emph{$d_k^-$-convex set} 
if it is isomorphic to the poset obtained by removing the top element from $d_k(1)$.

\begin{defi}
\label{def:dc_poset}
(See \cite{PS}, \cite{Okada})
A poset $P$ is \emph{$d$-complete} if it satisfies 
the following three conditions for every $k \ge 3$:
\begin{enumerate}
\item[(D1)]
If $I$ is a $d_k^-$-convex set, then there exists an element $u$ 
such that $u$ covers the maximal elements of $I$ and 
$I \cup \{ u \}$ is a $d_k$-interval.
\item[(D2)]
If $I = [v,u]$ is a $d_k$-interval and the top $u$ covers $u'$ in $P$, 
then $u' \in I$.
\item[(D3)]
There are no $d_k^-$-convex sets which differ only in the minimal elements.
\end{enumerate}
\end{defi}

It is clear that rooted trees, viewed as posets with their roots being the maximum elements, 
are $d$-complete posets.

\begin{exam}
\label{ex:shape}
For a partition $\lambda$, let $D(\lambda)$ be the Young diagram of $\lambda$ given by
$$
D(\lambda) = \{ (i,j) \in \Int^2 : i \ge 1, \ 1 \le j \le \lambda_i \}.
$$
For a strict partition $\mu$, let $S(\mu)$ be the shifted Young diagram of $\mu$ given by
$$
S(\mu) = \{ (i,j) \in \Int^2 : i \ge 1, \ i \le j \le \mu_i + i-1 \}.
$$
We endow $\Int^2$ with a poset structure by defining 
\begin{equation}
\label{eq:Z2-ordering}
\text{$(i,j) \ge (i',j')$ if $i \le i'$ and $j \le j'$.}
\end{equation}
Then we regard the Young diagram $D(\lambda)$ and the shifted Young diagram $S(\mu)$ 
as induced subposets of $\Int^2$.
The resulting posets are called a \emph{shape} and a \emph{shifted shape} respectively.
It can be shown that shapes and shifted shapes are $d$-complete posets.
Figure~\ref{fig:Hasse} illustrates the Hasse diagrams of 
$D(5,4,2,1)$ and $S(5,4,2,1)$.
\begin{figure}[!htb]
\centering
\subcaptionbox{$D(5,4,2,1)$}{
\setlength{\unitlength}{1.5pt}
\begin{picture}(80,70)(0,-10)
\put(5,15){\circle*{3}}
\put(15,25){\circle*{3}}
\put(25,15){\circle*{3}}
\put(25,35){\circle*{3}}
\put(35,25){\circle*{3}}
\put(35,45){\circle*{3}}
\put(45,15){\circle*{3}}
\put(45,35){\circle*{3}}
\put(55,5){\circle*{3}}
\put(55,25){\circle*{3}}
\put(65,15){\circle*{3}}
\put(75,5){\circle*{3}}
\put(5,15){\line(1,1){10}}
\put(15,25){\line(1,-1){10}}
\put(15,25){\line(1,1){10}}
\put(25,15){\line(1,1){10}}
\put(25,35){\line(1,-1){10}}
\put(25,35){\line(1,1){10}}
\put(35,25){\line(1,-1){10}}
\put(35,25){\line(1,1){10}}
\put(35,45){\line(1,-1){10}}
\put(45,15){\line(1,-1){10}}
\put(45,15){\line(1,1){10}}
\put(45,35){\line(1,-1){10}}
\put(55,5){\line(1,1){10}}
\put(55,25){\line(1,-1){10}}
\put(65,15){\line(1,-1){10}}
\end{picture}
}
\subcaptionbox{$S(5,4,2,1)$}{
\setlength{\unitlength}{1.5pt}
\begin{picture}(80,70)(-15,0)
\put(5,5){\circle*{3}}
\put(5,25){\circle*{3}}
\put(5,45){\circle*{3}}
\put(5,65){\circle*{3}}
\put(15,15){\circle*{3}}
\put(15,35){\circle*{3}}
\put(15,55){\circle*{3}}
\put(25,25){\circle*{3}}
\put(25,45){\circle*{3}}
\put(35,15){\circle*{3}}
\put(35,35){\circle*{3}}
\put(45,25){\circle*{3}}
\put(5,5){\line(1,1){10}}
\put(5,25){\line(1,-1){10}}
\put(5,25){\line(1,1){10}}
\put(5,45){\line(1,-1){10}}
\put(5,45){\line(1,1){10}}
\put(5,65){\line(1,-1){10}}
\put(15,15){\line(1,1){10}}
\put(15,35){\line(1,-1){10}}
\put(15,35){\line(1,1){10}}
\put(15,55){\line(1,-1){10}}
\put(25,25){\line(1,-1){10}}
\put(25,25){\line(1,1){10}}
\put(25,45){\line(1,-1){10}}
\put(35,15){\line(1,1){10}}
\put(35,35){\line(1,-1){10}}
\end{picture}
}
\caption{Shape and shifted shape}
\label{fig:Hasse}
\end{figure}
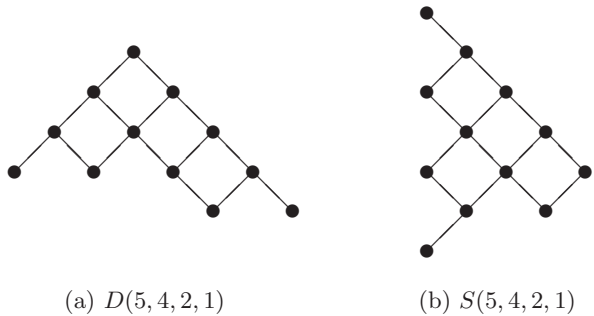
\end{exam}

\begin{exam}
\label{ex:swivel}
Let $P$ be a subset of $\Int^2$ given by
$$
P =
\left\{ 
 \begin{array}{l}
  (1,1), (1,2), (1,3), (1,4), (1,5), (2,3), (2,4), (2,5), \\
  (3,4), (3,5), (3,6), (4,4), (4,5), (4,6), (4,7), (4,8)
 \end{array}
\right\}.
$$
If we regard $P$ as an induced subposet of $\Int^2$ with ordering given by \eqref{eq:Z2-ordering}, 
then $P$ is a $d$-complete poset, called a \emph{swivel}.
See Figure~\ref{fig:Hasse_swivel} for the Hasse diagram of $P$.
\begin{figure}[!htb]
$$
\setlength{\unitlength}{1.5pt}
\begin{picture}(40,100)
\put(0,100){\circle*{3}}
\put(0,40){\circle*{3}}
\put(10,90){\circle*{3}}
\put(10,70){\circle*{3}}
\put(10,50){\circle*{3}}
\put(10,30){\circle*{3}}
\put(20,80){\circle*{3}}
\put(20,60){\circle*{3}}
\put(20,40){\circle*{3}}
\put(20,20){\circle*{3}}
\put(30,70){\circle*{3}}
\put(30,50){\circle*{3}}
\put(30,30){\circle*{3}}
\put(30,10){\circle*{3}}
\put(40,60){\circle*{3}}
\put(40,0){\circle*{3}}
\put(0,100){\line(1,-1){10}}
\put(0,40){\line(1,1){10}}
\put(0,40){\line(1,-1){10}}
\put(10,90){\line(1,-1){10}}
\put(10,70){\line(1,1){10}}
\put(10,70){\line(1,-1){10}}
\put(10,50){\line(1,1){10}}
\put(10,50){\line(1,-1){10}}
\put(10,30){\line(1,1){10}}
\put(10,30){\line(1,-1){10}}
\put(20,80){\line(1,-1){10}}
\put(20,60){\line(1,1){10}}
\put(20,60){\line(1,-1){10}}
\put(20,40){\line(1,1){10}}
\put(20,40){\line(1,-1){10}}
\put(20,20){\line(1,1){10}}
\put(20,20){\line(1,-1){10}}
\put(30,70){\line(1,-1){10}}
\put(30,50){\line(1,1){10}}
\put(30,10){\line(1,-1){10}}
\end{picture}
$$
\caption{Swivel}
\label{fig:Hasse_swivel}
\end{figure}
\end{exam}

A poset $P$ is called \emph{connected} if its Hasse diagram is a connected graph.
It is easy to see that, if $P$ is a $d$-complete poset, 
then each connected component of $P$ is $d$-complete.
(For our purpose to prove Theorem~\ref{thm:main}, 
there is no harm in assuming that a $d$-complete poset is connected.)

\begin{prop}
\label{prop:dc_connected}
\textup{(\cite[\S 3]{Proc1})}
Let $P$ be a $d$-complete poset.
If $P$ is connected, then $P$ has a unique maximal element.
\end{prop}

Let $P$ be a finite poset.
The \emph{top forest} $\Gamma$ of $P$ is a full subgraph of 
the Hasse diagram of $P$, whose vertex set consists of all elements $x \in P$ 
such that every $y \ge x$ is covered by at most one other element. 
Note that two elements $x$ and $y \in \Gamma$ are connected by an edge in $\Gamma$ 
if and only if $x$ covers $y$ or $y$ covers $x$.
Then $\Gamma$ becomes a forest as a graph (\ie a disjoint union of trees).
If $P$ is a connected $d$-complete poset, then $\Gamma$ is a tree called 
the \emph{top tree} of $P$.
We can regard the top forest $\Gamma$ as a Dynkin diagram of a Kac--Moody group.
For example, the top trees of Figure~\ref{fig:Hasse} (a) and (b) are
of type $A_8$ and $D_6$ respectively. 
If $P$ is the swivel in Figure~\ref{fig:Hasse_swivel}, the top tree is
the Dynkin diagram of type $E_6$.
(See also Subections~2.2 and 4.3.)

Remark~\ref{rem:heap1} below will indicate how to combine the results 
of \cite{Proc2} and \cite{Stem2} to obtain the following statement.

\begin{prop}
\label{prop:dc_color}
\textup{(\cite[Proposition~8.6]{Proc2}, \cite[Proposition~3.1]{Stem2})}
Let $P$ be a $d$-complete poset and $\Gamma$ its top forest.
Let $I$ be a set of colors whose cardinality is the same as $\Gamma$.
Then a bijective labeling $c : \Gamma \to I$ can be uniquely extended to a map 
$c : P \to I$ satisfying the following three conditions:
\begin{enumerate}
\item[(C1)]
If $x$ and $y$ are incomparable, then $c(x) \neq c(y)$.
\item[(C2)]
If an interval $[v,u]$ is a chain, then the colors $c(x)$ \textup{(}$x \in [v,u]$\textup{)} are 
distinct.
\item[(C3)]
If $[v,u]$ is a $d_k$-interval then $c(v) = c(u)$.
\end{enumerate}
Moreover this map $c$ satisfies
\begin{enumerate}
\item[(C4)]
If $x$ covers $y$, then the nodes labeled by $c(x)$ and $c(y)$ are adjacent in $\Gamma$.
\item[(C5)]
If $c(x) = c(y)$ or the nodes labeled by $c(x)$ and $c(y)$ are adjacent in $\Gamma$, 
then $x$ and $y$ are comparable.
\end{enumerate}
Such a map $c : P \to I$ is called a \emph{$d$-complete coloring}.
\end{prop}

Let $P$ be a $d$-complete poset and 
$c : P \to I$ a $d$-complete coloring.
Let $\vectz = (z_i)_{i \in I}$ be indeterminates.
Given an order filter $F$ of $P$, we regard $P \setminus F$ as an induced subposet.
For a $(P \setminus F)$-partition $\sigma \in \AP(P \setminus F)$, we put
$$
\vectz^\sigma = \prod_{v \in P \setminus F} z_{c(v)}^{\sigma(v)}.
$$
We are interested in the multivariate generating function
$$
\sum_{\sigma \in \AP(P \setminus F)} \vectz^\sigma
$$
of $(P \setminus F)$-partitions.
For a subset $D$ of $P$, we write
$$
\vectz[D] = \prod_{v \in D} z_{c(v)}.
$$
Instead of giving a definition of hooks $H_P(u) \subset P$ for a general 
$d$-complete poset $P$, we define associated monomials $\vectz[H_P(u)]$ 
directly by induction as follows:

\begin{defi}
\label{def:hook}
Let $P$ be a $d$-complete poset with $d$-complete coloring $c : P \to I$.
\begin{enumerate}
\item[(i)]
If $u$ is not the top of any $d_k$-interval, then we define
$$
\vectz[H_P(u)] = \prod_{w \le u} z_{c(w)}.
$$
\item[(ii)]
If $u$ is the top of a $d_k$-interval $[v,u]$, then we define
$$
\vectz[H_P(u)]
 =
\frac{ \vectz[H_P(x)] \cdot \vectz[H_P(y)] }
     { \vectz[H_P(v)] },
$$
where $x$ and $y$ are the sides of $[v,u]$.
\end{enumerate}
\end{defi}

\begin{exam}
\label{ex:shape_hc}
Let $P = D(\lambda)$ be the shape corresponding to a partition $\lambda$.
Then the top tree $\Gamma$ of $D(\lambda)$ is given by
$$
\Gamma = \{ (1,j) : 1 \le j \le \lambda_1 \} \cup \{ (i,1) : 1 \le i \le \lambda'_1 \},
$$
where $\lambda'_1$ is the number of cells in the first column of the Young diagram $D(\lambda)$.
A $d$-complete coloring $c: D(\lambda) \to I = \{ -(\lambda'_1-1), \dots, -1, 0, 1, \dots, \lambda_1-1 \}$ 
is given by
$$
c(i,j) = j-i.
$$
The classical definition of the hook $H_{D(\lambda)}(i,j) \subset D(\lambda)$ 
at $u = (i,j)$ in $D(\lambda)$ is as follows:
$$
H_{D(\lambda)}(i,j)
 =
\{ (i,j) \}
\cup
\{ (i,l) \in D(\lambda) : l > j \}
\cup
\{ (k,j) \in D(\lambda) : k > i \}.
$$
Then the hook monomial $\vectz[H_{D(\lambda)}(i,j)]$ in Definition~\ref{def:hook} 
and the hook $H_{D(\lambda)}(i,j)$ are related as 
$$
\vectz[H_{D(\lambda)}(i,j)]
 =
\prod_{(p,q)\in H_{D(\lambda)}(i,j) } z_{c(p,q)}.
$$
\end{exam}

\begin{exam}
\label{ex:shifted_hc}
Let $P = S(\mu)$ be the shifted shape corresponding to a strict partition $\mu$ of length $\ge 2$.
Then the top tree $\Gamma$ of $S(\mu)$ is given by
$$
\Gamma = \{ (1,j) : 1 \le j \le \mu_1 \} \cup \{ (2, 2) \},
$$
and a $d$-complete coloring $c : S(\mu) \to I = \{ 0, 0', 1, 2, \dots, \mu_1 - 1 \}$ is given by
$$
c(i,j)
 =
\begin{cases}
 j-i &\text{if $i < j$,} \\
 0 &\text{if $i=j$ and $i$ is odd,} \\
 0' &\text{if $i=j$ and $i$ is even.}
\end{cases}
$$
The (shifted) hook $H_{S(\mu)}(i,j)$ of $u = (i,j)$ in $S(\mu)$ is the subset of $S(\mu)$ defined by
\begin{align*}
H_{S(\mu)}(i,j)
 &=
\{ (i,j) \}
\cup
\{ (i,l) \in S(\mu) : l > j \}
\cup
\{ (k,j) \in S(\mu) : k > j \}
\\
&\quad
\cup
\{ (j+1,l) \in S(\mu) : l > j \}.
\end{align*}
For example, if $\mu = (5,4,2,1)$ and $(i,j) = (1,2)$, then the corresponding hook is given by
$$
H_{S(5,4,2,1)}(1,2)
 =
\{ (1,2), (1,3), (1,4), (1,5), (2,2), (3,3), (3,4) \},
$$
and the hook monomial is 
$\vectz[H_{S(5,4,2,1)}(1,2)]= z_{0'}z_{0} z_1^2 z_2 z_3 z_4$.
\end{exam}

\subsection{%
$d$-Complete posets and Weyl groups
}

Let $P$ be a connected $d$-complete poset with top tree $\Gamma$.
We regard $\Gamma$ as a (simply-laced) Dynkin diagram with node set $I$ 
and the $d$-complete coloring as a map $c : P \to I$.
Let $A = (a_{ij})_{i,j \in I}$ be the generalized Cartan matrix of $\Gamma$ given by
$$
a_{ij}
 =
\begin{cases}
 2 &\text{if $i=j$,} \\
 -1 &\text{if $i \neq j$ and $i$ and $j$ are adjacent in $\Gamma$,} \\
 0 &\text{otherwise.}
\end{cases}
$$
We fix the following data associated to $A$:
\begin{itemize}
\item
a free $\Int$-module $\Lambda$ of finite rank at least $\# \Gamma$, called the \emph{weight lattice},
\item
a linearly independent subset $\Pi = \{ \alpha_i : i \in I \}$ of $\Lambda$, called the \emph{simple roots},
\item
a subset $\Pi^\vee = \{ \alpha^\vee_i : i \in I \}$ of the dual lattice $\Lambda^* = \Hom_\Int(\Lambda,\Int)$, 
called the \emph{simple coroots},
\item
a subset $\{ \lambda_i : i \in I \}$ of $\Lambda$, called the \emph{fundamental weights},
\end{itemize}
such that
$$
\langle \alpha^\vee_i, \alpha_j \rangle = a_{ij},
\quad
\langle \alpha^\vee_i, \lambda_j \rangle = \delta_{ij},
$$
where $\langle \quad, \quad \rangle : \Lambda^* \times \Lambda \to \Int$ is the canonical pairing.
Let $W$ be the corresponding 
(Kac--Moody) Weyl group generated by the simple reflections $\{ s_i : i \in I \}$, 
where $s_i$ acts on $\Lambda$ and $\Lambda^*$ by the rule
$$
s_i(\lambda) = \lambda - \langle \alpha_i^\vee, \lambda \rangle \alpha_i
\quad(\lambda \in \Lambda),
\quad
s_i(\lambda^\vee) = \lambda^\vee - \langle \lambda^\vee, \alpha_i \rangle \alpha^\vee_i
\quad(\lambda^\vee \in \Lambda^*).
$$
Then $W$ is a Coxeter group, and we have the length function $l$ and the Bruhat order $<$ on $W$.
The set of real roots $\Phi$ and the set of real coroots $\Phi^\vee$ are defined by 
$\Phi = W \Pi$ and $\Phi^\vee = W \Pi^\vee$ respectively.
The set of simple roots $\Pi$ (\resp the set of simple coroots $\Pi^\vee$) determines 
the decomposition of $\Phi$ (\resp $\Phi^\vee$) into the positive system $\Phi_+$ (\resp $\Phi^\vee_+$) and 
the negative system $\Phi_-$ (\resp $\Phi^\vee_-$).
We introduce the standard partial ordering on $\Phi_+$ (\resp $\Phi^\vee_+$) by setting 
$\alpha > \beta$ if $\alpha - \beta$ is a sum of simple roots $\{ \alpha_i : i \in I \}$ 
(\resp $\alpha^\vee > \beta^\vee$ if $\alpha^\vee - \beta^\vee$ is a sum of simple coroots 
$\{ \alpha^\vee_i : i \in I \}$).

For $p \in P$, we put
$$
\alpha(p) = \alpha_{c(p)},
\quad
\alpha^\vee(p) = \alpha^\vee_{c(p)},
\quad
s(p) = s_{c(p)}.
$$
Let $\alpha_P$ and $\lambda_P$ be the simple root and the fundamental weight 
corresponding to the color $i_P$ of the maximum element of $P$. 

Take a linear extension and label the elements of $P$ with $p_1, \cdots, p_N$ ($N= \# P$) 
so that $p_i < p_j$ in $P$ implies $i<j$.
Then we construct an element $w_P \in W$ by putting
$$
w_P = s(p_1) s(p_2) \cdots s(p_N).
$$
A Weyl group element $w \in W$ is called \emph{$\lambda$-minuscule} 
if there exists a reduced expression $w = s_{i_1} \cdots s_{i_l}$ such that
$$
\langle \alpha^\vee_{i_k}, s_{i_{k+1}} \cdots s_{i_l} \lambda \rangle = 1
\quad(1 \le k \le l),
$$
or equivalently
$$
s_{i_k} \cdots s_{i_l} \lambda = \lambda - \alpha_{i_k} - \cdots - \alpha_{i_l}.
$$
A element $w \in W$ is called \emph{fully commutative} 
if any reduced expression of $w$ can be obtained from any other by using only 
the Coxeter relations of the form $st = ts$.

\begin{prop}
\label{prop:w_P}
\textup{(See \cite{Proc2} and \cite[Proposition~2.1]{Stem2})}
Let $P$ be a connected $d$-complete poset.
Then the Weyl group element $w_P \in W$ is $\lambda_P$-minuscule and hence fully commutative.
\end{prop}

If $p = p_k \in P$, then we define 
\begin{align*}
\beta(p_k) &= s(p_1) \cdots s(p_{k-1}) \alpha(p_k),
\\
\gamma(p_k) &= s(p_N) \cdots s(p_{k+1}) \alpha(p_k),
\\
\gamma^\vee(p_k) &= s(p_N) \cdots s(p_{k+1}) \alpha^\vee(p_k).
\end{align*}
It follows from Proposition~\ref{prop:w_P} that, 
for each $p \in P$, the roots $\beta(p)$, $\gamma(p)$ and the coroot $\gamma^\vee(p)$ 
are independent of the choices of linear extensions.
For a Weyl group element $w \in W$, we put
$$
\Phi(w) = \Phi_+ \cap w \Phi_-,
\quad
\Phi^\vee(w) = \Phi^\vee_+ \cap w \Phi^\vee_-.
$$
Then it is well-known (see \cite[\S5.6]{H}) that
$$
\Phi(w_P) = \{ \beta(p) : p \in P \},
\quad
\Phi^\vee(w_P^{-1}) = \{ \gamma^\vee(p) : p \in P \}.
$$
Moreover we have

\begin{prop}
\label{prop:root}
Let $P$ be a connected $d$-complete poset. Then we have
\begin{enumerate}
\item[(a)]
\textup{(See \cite[Proposition~3.1 and Theorem~5.5]{Stem2})}
The poset $P$ is isomorphic to the order dual of $\Phi^\vee(w_P^{-1})$ 
with the standard coroot ordering on $\Phi^\vee_+$.
\item[(b)]
\textup{(See \cite[Lemma~IV]{Proc3})}
Under the identification $z_i = e^{\alpha_i}$ \textup{(}$i \in I$\textup{)}, we have
$$
\vectz[H_P(p)] = e^{\beta(p)}
\quad(p \in P).
$$
\item[(c)]
\textup{(See \cite[Proposition~5.1]{Stem2})}
We have
$$
\langle \gamma^\vee(p), \lambda_P \rangle = 1
\quad(p \in P).
$$
\end{enumerate}
\end{prop}

Let $W_{\lambda_P}$ be the stabilizer of $\lambda_P$ in $W$.
Then $W_{\lambda_P}$ is the maximal parabolic subgroup corresponding to $I \setminus \{ i_P \}$.
Let $W^{\lambda_P}$ be the set of minimum length coset representatives of $W/W_{\lambda_P}$.
For a subset $D = \{ p_{i_1}, \cdots, p_{i_r} \}$ ($i_1 < \cdots < i_r$) of $P$, we define
\begin{equation}
\label{eq:w_D}
w_D = s(p_{i_1}) \cdots s(p_{i_r}).
\end{equation}
Since $w_P$ is fully commutative (Proposition~\ref{prop:w_P}), 
we see that $w_D$ is independent of the choices of linear extensions of $P$.

\begin{prop}
\label{prop:Weyl}
Let $P$ be a connected $d$-complete poset.
Then we have
\begin{enumerate}
\item[(a)]
\textup{(See \cite[Proposition~I]{Proc3})}
The map $F \mapsto w_F$ gives a poset isomorphism 
from the set of all order filters of $P$ ordered by inclusion 
to the Bruhat interval $[e,w_P]$ in $W^{\lambda_P}$.
\item[(b)]
\textup{(See \cite[Remark~2.7 (b)]{Stem2})}
If $F$ is an order filter of $P$, then $w_F$ is $\lambda_P$-minuscule, and
$$
w_F \lambda_P = \lambda_P - \sum_{p \in F} \alpha(p).
$$
\end{enumerate}
\end{prop}

\begin{rema}
\label{rem:heap1}
Let $A = (a_{ij})_{i,j \in I}$ be a symmetrizable generalized Cartan matrix, 
and $\Gamma$ the corresponding Dynkin diagram with node set $I$.
Let $W$ be the associated Kac--Moody Weyl group generated by $\{ s_i : i \in I \}$.
Given a (not necessarily reduced) expression $s_{i_1} s_{i_2} \dots s_{i_N}$ of an element $w \in W$ 
in simple reflections, we can define a poset $H$, called the \emph{heap}, as follows 
(see \cite{Stem1}).
The poset $H$ consists of the ground set $\{ 1, 2, \dots, N \}$ and the partial ordering 
$\preceq$ obtained by taking the transitive closure of the relations given by
$$
\text{$a \prec b$ if $a < b$ and either $s_{i_a} s_{i_b} \neq s_{i_b} s_{i_a}$ or $i_a = i_b$.}
$$
The heap $H$ has a natural labeling (coloring) $c: H \to I$ given by $c(a) = i_a$.
If $w \in W$ is fully commutative, then the heap defined by a reduced expression of $w$ 
is independent of the choices of reduced expressions.
In this case we denote the resulting heap by $H(w)$.

If $a_{ij} = a_{ji} \in \{ 0, -1 \}$ for any pair of vertices $i$, $j \in I$ ($i \neq j$), 
then we say that $A$, $W$ and $H(w)$ are \emph{simply-laced}; 
otherwise they are \emph{multiply-laced}.
For a dominant weight $\lambda$ and a $\lambda$-minuscule element $w$ in the simply-laced Weyl group $W$,
Proctor \cite{Proc2} proved that the interval $[e,w]$ of $W^\lambda$ with respect to 
the Bruhat order is a distributive lattice, 
and that the order dual of the induced subposet consisting of all join-irreducible elements of $[e,w]$ 
is a $d$-complete poset colored by $I$.
Stembridge \cite{Stem2} extends Proctor's result to any symmetrizable Kac--Moody Weyl group 
in the setting of heaps.
The results of \cite{Proc2} and \cite{Stem2} may be combined by using 
\cite[Theorem~3.2 (c)]{Stem1} to produce Proposition~\ref{prop:dc_color} above.
Henceforth the simply-laced case will be referred to with the ``$d$-complete poset'' terminology 
and the general case will be referred to with the ``heap'' terminology.
Every $d$-complete poset $P$ is isomorphic to the simply-laced heap $H(w_P)$. 
In general, if $w \in W$ is dominant minuscule, \ie $\lambda$-minuscule for some dominant weight $\lambda$, 
the corresponding heap $H(w)$ is isomorphic (as a unlabeled poset) to a $d$-complete poset.
See \cite[Sections~3 and 4]{Stem2}.

The propositions in this section hold literally for heaps $H(w)$ of dominant minuscule elements, 
except for Proposition~\ref{prop:dc_color} and Proposition~\ref{prop:root} (b).
The latter half of Proposition~\ref{prop:dc_color} holds for heaps, 
\ie the labeling $c : H(w) \to I$ satisfies (C4) and (C5).
And we adopt Proposition~\ref{prop:root} (b) as a definition of the hook monomial for $H(w)$.
Below is an example of a non-simply-laced heap.
\end{rema}

\begin{exam} (Non-simply-laced heap)
\label{ex:shifted-B}
Let $\mu = (\mu_1, \dots, \mu_l)$ be a strict partition of length $l$.
Then the shifted shape $S(\mu)$ can be regarded as 
the heap associated to a dominant minuscule element of the Weyl group of type $B$,
which is not simply-laced.
Put $m = \mu_1$ and let $W$ be the Weyl group generated by $s_0, s_1, \dots, s_{m-1}$ 
subject to the relations
\begin{gather*}
s_i^2 = 1 \quad(i=0, 1, \dots, m-1),
\\
s_0 s_1 s_0 s_1 = s_1 s_0 s_1 s_0,
\\
s_i s_{i+1} s_i = s_{i+1} s_i s_{i+1} \quad(i=1, 2, \dots, m-2),
\\
s_i s_j = s_j s_i \quad(|i-j| \ge 2).
\end{gather*}
Then we define an element $w_\mu \in W$ by putting
$$
w_\mu = (s_{\mu_l-1} \cdots s_1 s_0) \cdots (s_{\mu_2-1} \cdots s_1 s_0) (s_{\mu_1-1} \cdots s_1 s_0).
$$
Then it can be shown that $w_\mu$ is $\lambda_0$-minuscule, where $\lambda_0$ is the fundamental weight 
corresponding to $s_0$, and that the map
\begin{equation}
\label{eq:S=H}
S(\mu) \ni (i,j) \mapsto \mu_l + \dots + \mu_i - j + i \in H(w_\mu)
\end{equation}
gives a poset isomorphism.
We identify the ground set of $H(w_\mu)$ with $S(\mu)$ via the isomorphism \eqref{eq:S=H}.
Then the natural labeling $c' : S(\mu) \to \{ 0, 1, \dots, m-1 \}$ of $H(w_\mu)$ is given by 
$c'(i,j) = j-i$, which is different from the $d$-complete coloring of $S(\mu)$ given 
in Example~\ref{ex:shifted_hc}.
And the ``hook'' $H'_{S(\mu)}(v)$ (see \cite[Definition~4.8]{N2}) is defined by
\begin{align*}
H'_{S(\mu)}(i,j)
 &=
\{ (i,j) \}
 \cup
\{ (i,l) \in S(\mu) : l > j \}
 \cup
\{ (k,j) \in S(\mu) : k > i \}
\\
 &\quad
\cup
\begin{cases}
 \{ (i,i) \} \cup \{ (j,l) \in S(\mu) : l \ge j \} &\text{if $i<j$ and $(j,j) \in S(\mu)$,} \\
 \emptyset &\text{otherwise,}
\end{cases}
\end{align*}
which is also different from the shifted hook given in Example~\ref{ex:shifted_hc}.
(See also Example~\ref{ex:heap}.)
\end{exam}
\section{%
Excited diagrams
}

In this section we introduce the notion of excited and $K$-theoretical excited diagrams 
in a $d$-complete poset and study their properties.

\subsection{%
Excited diagrams
}

First we generalize the notion of exited diagrams for Young diagram and shifted Young diagram, 
which were introduced by Ikeda--Naruse \cite{IN1} and Kreiman \cite{Kr1}, \cite{Kr2} independently, 
to a general $d$-complete posets.
And we give a generalization of backward movable positions or excited peaks 
introduced in \cite{IN2}, \cite{KN} and \cite{MPP}.

Let $P$ be a $d$-complete poset with top forest $\Gamma$ and $d$-complete coloring $c : P \to I$.
For a subset $D \subset P$ and a color $i \in I$, we put
$$
D_i = \{ x \in D : c(x) = i \}.
$$
For $i \in I$, let $N_i$ be the subset of $P$ consisting of element $x \in P$ whose color 
$c(x)$ is adjacent to $i$ in the Dynkin diagram $\Gamma$.
Note that, if $[v,u]$ is a $d_k$-interval, then $[v,u] \cap N_{c(u)}$ consists of elements $x \in [v,u]$ 
such that $x$ is covered by $u$ or covers $v$.

\begin{defi}
\label{def:excited}
Let $P$ be a $d$-complete poset and let $F$ be an order filter of $P$.
\begin{enumerate}
\item[(a)]
Let $D$ be a subset of $P$ and $u \in D$.
We say that $u$ is \emph{$D$-active} 
if there exists an element $v \in (P \setminus D)_{c(u)}$ 
such that $v < u$, $[v,u]$ is a $d_k$-interval and
$$
[v,u] \cap D \cap N_{c(u)} = \emptyset.
$$
\item[(b)]
Let $D$ be a subset of $P$ and $u \in D$.
If $u$ is $D$-active, then we define $\alpha_u(D)$ to be the subset of $P$ obtained from $D$ by 
replacing $u \in D$ by the bottom element $v$ of the $d_k$-interval $[v,u]$.
We call this replacement an (ordinary) \emph{elementary excitation}.
\item[(c)]
An \emph{excited diagram} of $F$ in $P$ is a subset of $P$ 
obtained from $F$ after a sequence of elementary excitations on active elements. 
Let $\EE_P(F)$ be the set of all excited diagrams of $F$ in $P$.
\item[(d)]
To an excited diagram $D \in \EE_P(F)$ we associate a subset $B(D) \subset P$ as follows:
If $D = F$, then $B(F) = \emptyset$.
If $D$ is an excited diagram with an active element $u$, then we define
$$
B(\alpha_u(D))
 = 
\left( B(D) \setminus ([v,u] \cap N_{c(u)}) \right) \cup \{ u \},
$$
where $[v,u]$ is the $d_k$-interval with top element $u$.
We call $B(D)$ the set of \emph{excited peaks} of $D$.
(We will show that $B(D)$ is a well-defined subset of $P \setminus D$ in Proposition~\ref{prop:BD}.)
\end{enumerate}
\end{defi}

In general, if two elements $u \in D$ and $v \not\in D$ with $v < u$ 
have the same color $i = c(v) = c(u)$ and satisfy $[v,u] \cap D \cap N_i = \emptyset$, 
then $D \setminus \{ u \} \cup \{ v \}$ is obtained from $D$ by a sequence 
of elementary excitations.

If $P$ is a shape or a shifted shape, our definition above coincides with the definitions of 
elementary excitations in \cite{IN1}, \cite{IN2}, 
or ladder moves in \cite{Kr1}, \cite{Kr2}, 
and backward movable positions in \cite{IN2}, \cite{KN} 
or excited peaks in \cite{MPP} (only for a shape).

\begin{exam}
\label{ex:shape_excited}
If $P = D(5,4,2,1)$ is the shape corresponding to a partition $(5,4,2,1)$ and 
$F = D(3,1)$, then there are $7$ excited diagrams in $\EE_P(F)$ shown in Figure~\ref{fig:shape_excited}.
In Figure~\ref{fig:shape_excited} (and Figures~\ref{fig:swivel_excited}) 
the shaded cells form an exited diagram and a cell with $\times$ is an excited peak.
And the arrow $D \longrightarrow D'$ means that $D'$ is obtained from $D$ by an elementary excitation.
\begin{figure}[t]
$$
\setlength{\unitlength}{1.1pt}
\begin{CD}
\raisebox{-20pt}{
\begin{picture}(50,40)
\fboxsep=0mm
\put(0,30){\colorbox[gray]{0.7}{\makebox(10,10){}}}
\put(10,30){\colorbox[gray]{0.7}{\makebox(10,10){}}}
\put(20,30){\colorbox[gray]{0.7}{\makebox(10,10){}}}
\put(0,20){\colorbox[gray]{0.7}{\makebox(10,10){}}}
\put(0,40){\line(1,0){50}}
\put(0,30){\line(1,0){50}}
\put(0,20){\line(1,0){40}}
\put(0,10){\line(1,0){20}}
\put(0,0){\line(1,0){10}}
\put(0,0){\line(0,1){40}}
\put(10,0){\line(0,1){40}}
\put(20,10){\line(0,1){30}}
\put(30,20){\line(0,1){20}}
\put(40,20){\line(0,1){20}}
\put(50,30){\line(0,1){10}}
\end{picture}
}
@>>>
\raisebox{-20pt}{
\begin{picture}(50,40)
\fboxsep=0mm
\put(0,30){\colorbox[gray]{0.7}{\makebox(10,10){}}}
\put(10,30){\colorbox[gray]{0.7}{\makebox(10,10){}}}
\put(30,20){\colorbox[gray]{0.7}{\makebox(10,10){}}}
\put(0,20){\colorbox[gray]{0.7}{\makebox(10,10){}}}
\put(0,40){\line(1,0){50}}
\put(0,30){\line(1,0){50}}
\put(0,20){\line(1,0){40}}
\put(0,10){\line(1,0){20}}
\put(0,0){\line(1,0){10}}
\put(0,0){\line(0,1){40}}
\put(10,0){\line(0,1){40}}
\put(20,10){\line(0,1){30}}
\put(30,20){\line(0,1){20}}
\put(40,20){\line(0,1){20}}
\put(50,30){\line(0,1){10}}
\put(20,30){\makebox(10,10){$\times$}}
\end{picture}
}
@>>>
\raisebox{-20pt}{
\begin{picture}(50,40)
\fboxsep=0mm
\put(0,30){\colorbox[gray]{0.7}{\makebox(10,10){}}}
\put(20,20){\colorbox[gray]{0.7}{\makebox(10,10){}}}
\put(30,20){\colorbox[gray]{0.7}{\makebox(10,10){}}}
\put(0,20){\colorbox[gray]{0.7}{\makebox(10,10){}}}
\put(0,40){\line(1,0){50}}
\put(0,30){\line(1,0){50}}
\put(0,20){\line(1,0){40}}
\put(0,10){\line(1,0){20}}
\put(0,0){\line(1,0){10}}
\put(0,0){\line(0,1){40}}
\put(10,0){\line(0,1){40}}
\put(20,10){\line(0,1){30}}
\put(30,20){\line(0,1){20}}
\put(40,20){\line(0,1){20}}
\put(50,30){\line(0,1){10}}
\put(10,30){\makebox(10,10){$\times$}}
\end{picture}
}
\\[5pt]
@VVV
@VVV
@VVV
\\[5pt]
\raisebox{-20pt}{
\begin{picture}(50,40)
\fboxsep=0mm
\put(0,30){\colorbox[gray]{0.7}{\makebox(10,10){}}}
\put(10,30){\colorbox[gray]{0.7}{\makebox(10,10){}}}
\put(20,30){\colorbox[gray]{0.7}{\makebox(10,10){}}}
\put(10,10){\colorbox[gray]{0.7}{\makebox(10,10){}}}
\put(0,40){\line(1,0){50}}
\put(0,30){\line(1,0){50}}
\put(0,20){\line(1,0){40}}
\put(0,10){\line(1,0){20}}
\put(0,0){\line(1,0){10}}
\put(0,0){\line(0,1){40}}
\put(10,0){\line(0,1){40}}
\put(20,10){\line(0,1){30}}
\put(30,20){\line(0,1){20}}
\put(40,20){\line(0,1){20}}
\put(50,30){\line(0,1){10}}
\put(0,20){\makebox(10,10){$\times$}}
\end{picture}
}
@>>>
\raisebox{-20pt}{
\begin{picture}(50,40)
\fboxsep=0mm
\put(0,30){\colorbox[gray]{0.7}{\makebox(10,10){}}}
\put(10,30){\colorbox[gray]{0.7}{\makebox(10,10){}}}
\put(30,20){\colorbox[gray]{0.7}{\makebox(10,10){}}}
\put(10,10){\colorbox[gray]{0.7}{\makebox(10,10){}}}
\put(0,40){\line(1,0){50}}
\put(0,30){\line(1,0){50}}
\put(0,20){\line(1,0){40}}
\put(0,10){\line(1,0){20}}
\put(0,0){\line(1,0){10}}
\put(0,0){\line(0,1){40}}
\put(10,0){\line(0,1){40}}
\put(20,10){\line(0,1){30}}
\put(30,20){\line(0,1){20}}
\put(40,20){\line(0,1){20}}
\put(50,30){\line(0,1){10}}
\put(0,20){\makebox(10,10){$\times$}}
\put(20,30){\makebox(10,10){$\times$}}
\end{picture}
}
@>>>
\raisebox{-20pt}{
\begin{picture}(50,40)
\fboxsep=0mm
\put(0,30){\colorbox[gray]{0.7}{\makebox(10,10){}}}
\put(20,20){\colorbox[gray]{0.7}{\makebox(10,10){}}}
\put(30,20){\colorbox[gray]{0.7}{\makebox(10,10){}}}
\put(10,10){\colorbox[gray]{0.7}{\makebox(10,10){}}}
\put(0,40){\line(1,0){50}}
\put(0,30){\line(1,0){50}}
\put(0,20){\line(1,0){40}}
\put(0,10){\line(1,0){20}}
\put(0,0){\line(1,0){10}}
\put(0,0){\line(0,1){40}}
\put(10,0){\line(0,1){40}}
\put(20,10){\line(0,1){30}}
\put(30,20){\line(0,1){20}}
\put(40,20){\line(0,1){20}}
\put(50,30){\line(0,1){10}}
\put(0,20){\makebox(10,10){$\times$}}
\put(10,30){\makebox(10,10){$\times$}}
\end{picture}
}
@>>>
\raisebox{-20pt}{
\begin{picture}(50,40)
\fboxsep=0mm
\put(10,20){\colorbox[gray]{0.7}{\makebox(10,10){}}}
\put(20,20){\colorbox[gray]{0.7}{\makebox(10,10){}}}
\put(30,20){\colorbox[gray]{0.7}{\makebox(10,10){}}}
\put(10,10){\colorbox[gray]{0.7}{\makebox(10,10){}}}
\put(0,40){\line(1,0){50}}
\put(0,30){\line(1,0){50}}
\put(0,20){\line(1,0){40}}
\put(0,10){\line(1,0){20}}
\put(0,0){\line(1,0){10}}
\put(0,0){\line(0,1){40}}
\put(10,0){\line(0,1){40}}
\put(20,10){\line(0,1){30}}
\put(30,20){\line(0,1){20}}
\put(40,20){\line(0,1){20}}
\put(50,30){\line(0,1){10}}
\put(0,30){\makebox(10,10){$\times$}}
\end{picture}
}
\end{CD}
$$
\caption{Excited diagrams of $D(3,1)$ in $D(5,4,2,1)$}
\label{fig:shape_excited}
\end{figure}
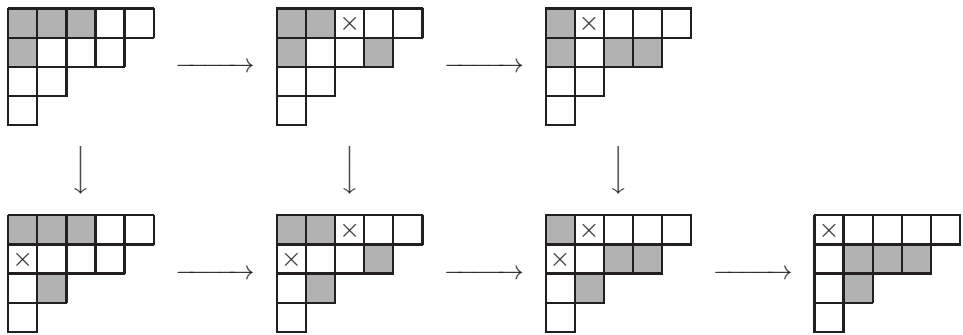
\end{exam}

\begin{exam}
\label{ex:swivel_excited}
If $P$ is the swivel given in Example~\ref{ex:swivel} and $F$ is the order filter consisting of three elements, 
then there are $8$ excited diagrams in $\EE_P(F)$ shown in Figure~\ref{fig:swivel_excited}.
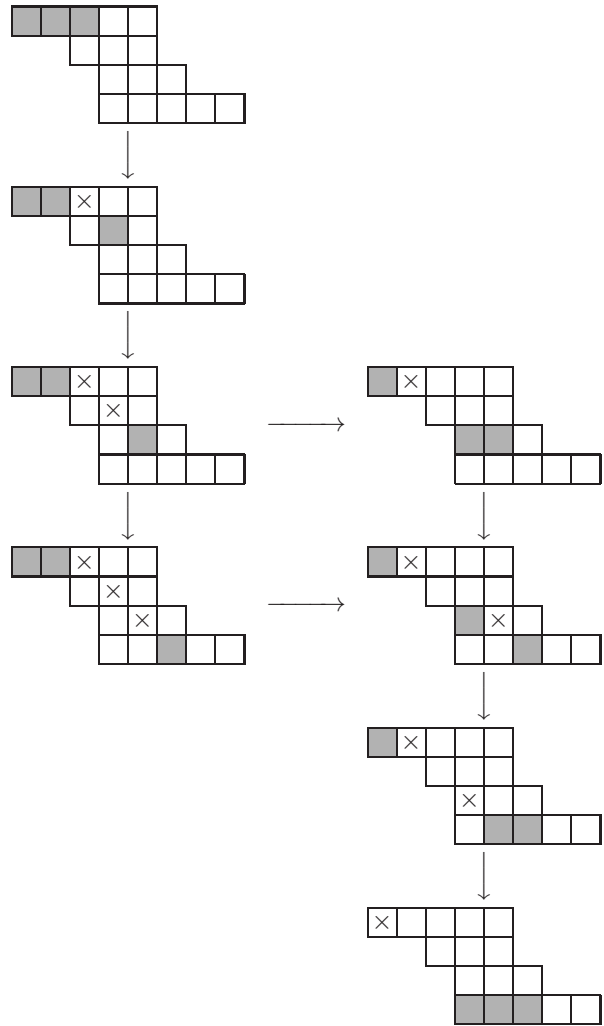
\begin{figure}[!htb]
$$
\setlength{\unitlength}{1.1pt}
\begin{CD}
\raisebox{-20pt}{
\begin{picture}(80,40)
\fboxsep=0mm
\put(0,30){\colorbox[gray]{0.7}{\makebox(10,10){}}}
\put(10,30){\colorbox[gray]{0.7}{\makebox(10,10){}}}
\put(20,30){\colorbox[gray]{0.7}{\makebox(10,10){}}}
\put(0,40){\line(1,0){50}}
\put(0,30){\line(1,0){50}}
\put(20,20){\line(1,0){40}}
\put(30,10){\line(1,0){50}}
\put(30,0){\line(1,0){50}}
\put(0,30){\line(0,1){10}}
\put(10,30){\line(0,1){10}}
\put(20,20){\line(0,1){20}}
\put(30,0){\line(0,1){40}}
\put(40,0){\line(0,1){40}}
\put(50,0){\line(0,1){40}}
\put(60,0){\line(0,1){20}}
\put(70,0){\line(0,1){10}}
\put(80,0){\line(0,1){10}}
\end{picture}
}
\\
@VVV
\\
\raisebox{-20pt}{
\begin{picture}(80,40)
\fboxsep=0mm
\put(0,30){\colorbox[gray]{0.7}{\makebox(10,10){}}}
\put(10,30){\colorbox[gray]{0.7}{\makebox(10,10){}}}
\put(30,20){\colorbox[gray]{0.7}{\makebox(10,10){}}}
\put(0,40){\line(1,0){50}}
\put(0,30){\line(1,0){50}}
\put(20,20){\line(1,0){40}}
\put(30,10){\line(1,0){50}}
\put(30,0){\line(1,0){50}}
\put(0,30){\line(0,1){10}}
\put(10,30){\line(0,1){10}}
\put(20,20){\line(0,1){20}}
\put(30,0){\line(0,1){40}}
\put(40,0){\line(0,1){40}}
\put(50,0){\line(0,1){40}}
\put(60,0){\line(0,1){20}}
\put(70,0){\line(0,1){10}}
\put(80,0){\line(0,1){10}}
\put(20,30){\makebox(10,10){$\times$}}
\end{picture}
}
\\
@VVV
\\
\raisebox{-20pt}{
\begin{picture}(80,40)
\fboxsep=0mm
\put(0,30){\colorbox[gray]{0.7}{\makebox(10,10){}}}
\put(10,30){\colorbox[gray]{0.7}{\makebox(10,10){}}}
\put(40,10){\colorbox[gray]{0.7}{\makebox(10,10){}}}
\put(0,40){\line(1,0){50}}
\put(0,30){\line(1,0){50}}
\put(20,20){\line(1,0){40}}
\put(30,10){\line(1,0){50}}
\put(30,0){\line(1,0){50}}
\put(0,30){\line(0,1){10}}
\put(10,30){\line(0,1){10}}
\put(20,20){\line(0,1){20}}
\put(30,0){\line(0,1){40}}
\put(40,0){\line(0,1){40}}
\put(50,0){\line(0,1){40}}
\put(60,0){\line(0,1){20}}
\put(70,0){\line(0,1){10}}
\put(80,0){\line(0,1){10}}
\put(20,30){\makebox(10,10){$\times$}}
\put(30,20){\makebox(10,10){$\times$}}
\end{picture}
}
@>>>
\raisebox{-20pt}{
\begin{picture}(80,40)
\fboxsep=0mm
\put(0,30){\colorbox[gray]{0.7}{\makebox(10,10){}}}
\put(30,10){\colorbox[gray]{0.7}{\makebox(10,10){}}}
\put(40,10){\colorbox[gray]{0.7}{\makebox(10,10){}}}
\put(0,40){\line(1,0){50}}
\put(0,30){\line(1,0){50}}
\put(20,20){\line(1,0){40}}
\put(30,10){\line(1,0){50}}
\put(30,0){\line(1,0){50}}
\put(0,30){\line(0,1){10}}
\put(10,30){\line(0,1){10}}
\put(20,20){\line(0,1){20}}
\put(30,0){\line(0,1){40}}
\put(40,0){\line(0,1){40}}
\put(50,0){\line(0,1){40}}
\put(60,0){\line(0,1){20}}
\put(70,0){\line(0,1){10}}
\put(80,0){\line(0,1){10}}
\put(10,30){\makebox(10,10){$\times$}}
\end{picture}
}
\\
@VVV @VVV
\\
\raisebox{-20pt}{
\begin{picture}(80,40)
\fboxsep=0mm
\put(0,30){\colorbox[gray]{0.7}{\makebox(10,10){}}}
\put(10,30){\colorbox[gray]{0.7}{\makebox(10,10){}}}
\put(50,0){\colorbox[gray]{0.7}{\makebox(10,10){}}}
\put(0,40){\line(1,0){50}}
\put(0,30){\line(1,0){50}}
\put(20,20){\line(1,0){40}}
\put(30,10){\line(1,0){50}}
\put(30,0){\line(1,0){50}}
\put(0,30){\line(0,1){10}}
\put(10,30){\line(0,1){10}}
\put(20,20){\line(0,1){20}}
\put(30,0){\line(0,1){40}}
\put(40,0){\line(0,1){40}}
\put(50,0){\line(0,1){40}}
\put(60,0){\line(0,1){20}}
\put(70,0){\line(0,1){10}}
\put(80,0){\line(0,1){10}}
\put(20,30){\makebox(10,10){$\times$}}
\put(30,20){\makebox(10,10){$\times$}}
\put(40,10){\makebox(10,10){$\times$}}
\end{picture}
}
@>>>
\raisebox{-20pt}{
\begin{picture}(80,40)
\fboxsep=0mm
\put(0,30){\colorbox[gray]{0.7}{\makebox(10,10){}}}
\put(30,10){\colorbox[gray]{0.7}{\makebox(10,10){}}}
\put(50,0){\colorbox[gray]{0.7}{\makebox(10,10){}}}
\put(0,40){\line(1,0){50}}
\put(0,30){\line(1,0){50}}
\put(20,20){\line(1,0){40}}
\put(30,10){\line(1,0){50}}
\put(30,0){\line(1,0){50}}
\put(0,30){\line(0,1){10}}
\put(10,30){\line(0,1){10}}
\put(20,20){\line(0,1){20}}
\put(30,0){\line(0,1){40}}
\put(40,0){\line(0,1){40}}
\put(50,0){\line(0,1){40}}
\put(60,0){\line(0,1){20}}
\put(70,0){\line(0,1){10}}
\put(80,0){\line(0,1){10}}
\put(10,30){\makebox(10,10){$\times$}}
\put(40,10){\makebox(10,10){$\times$}}
\end{picture}
}
\\
@. @VVV
\\
@.
\raisebox{-20pt}{
\begin{picture}(80,40)
\fboxsep=0mm
\put(0,30){\colorbox[gray]{0.7}{\makebox(10,10){}}}
\put(40,0){\colorbox[gray]{0.7}{\makebox(10,10){}}}
\put(50,0){\colorbox[gray]{0.7}{\makebox(10,10){}}}
\put(0,40){\line(1,0){50}}
\put(0,30){\line(1,0){50}}
\put(20,20){\line(1,0){40}}
\put(30,10){\line(1,0){50}}
\put(30,0){\line(1,0){50}}
\put(0,30){\line(0,1){10}}
\put(10,30){\line(0,1){10}}
\put(20,20){\line(0,1){20}}
\put(30,0){\line(0,1){40}}
\put(40,0){\line(0,1){40}}
\put(50,0){\line(0,1){40}}
\put(60,0){\line(0,1){20}}
\put(70,0){\line(0,1){10}}
\put(80,0){\line(0,1){10}}
\put(10,30){\makebox(10,10){$\times$}}
\put(30,10){\makebox(10,10){$\times$}}
\end{picture}
}
\\
@. @VVV
\\
@.
\raisebox{-20pt}{
\begin{picture}(80,40)
\fboxsep=0mm
\put(30,0){\colorbox[gray]{0.7}{\makebox(10,10){}}}
\put(40,0){\colorbox[gray]{0.7}{\makebox(10,10){}}}
\put(50,0){\colorbox[gray]{0.7}{\makebox(10,10){}}}
\put(0,40){\line(1,0){50}}
\put(0,30){\line(1,0){50}}
\put(20,20){\line(1,0){40}}
\put(30,10){\line(1,0){50}}
\put(30,0){\line(1,0){50}}
\put(0,30){\line(0,1){10}}
\put(10,30){\line(0,1){10}}
\put(20,20){\line(0,1){20}}
\put(30,0){\line(0,1){40}}
\put(40,0){\line(0,1){40}}
\put(50,0){\line(0,1){40}}
\put(60,0){\line(0,1){20}}
\put(70,0){\line(0,1){10}}
\put(80,0){\line(0,1){10}}
\put(0,30){\makebox(10,10){$\times$}}
\end{picture}
}
\end{CD}
$$
\caption{Excited diagrams in a swivel}
\label{fig:swivel_excited}
\end{figure}
\end{exam}

Let $P$ be a connected $d$-complete poset 
and $W$ the Weyl group corresponding to the top tree $\Gamma$ viewed as a Dynkin diagram. 
Fix a linear extension of $P$, 
\ie a labeling of the elements of $P$ with $p_1, \dots, p_N$ so that $p_i < p_j$ in $P$ 
implies $i < j$.
Then we can associate to a subset $D$ of $P$ a well-defined element $w_D \in W$ as in 
\eqref{eq:w_D}.
The following proposition gives a characterization of excited diagrams.

\begin{prop}
\label{prop:E=R}
Let $P$ be a connected $d$-complete poset and let $F$ be an order filter of $P$.
Then a subset $D \subset P$ is an excited diagram of $F$ in $P$ 
if and only if $\# D = \# F$ and $w_D = w_F$.
\end{prop}

To prove this proposition, we prepare two lemmas.

\begin{lemm}
\label{lem:comm}
Let $D$ be a subset of $P$ and $u$ and $v$ elements of $P$ such that $v < u$ and $c(u) = c(v) = i$.
Suppose $v = p_k$ and $u = p_l$.
If $[v,u] \cap D \cap N_i = \emptyset$, 
then we have $s(p_j) s(p_l) = s(p_l) s(p_j)$ for any $p_j \in D$ with $k < j < l$ 
\end{lemm}

\begin{proof}
follows from Property (C5) in Proposition~\ref{prop:dc_color}.
\end{proof}

\begin{lemm}
\label{lem:fc}
\textup{(See \cite[Section~3]{GK})}
Let $v \in W$ be a fully commutative element 
and $v = s_{i_1} \cdots s_{i_r} = s_{j_1} \cdots s_{j_r}$ be its reduced expressions.
Then we have
\begin{enumerate}
\item[(a)]
Let $i$, $j$ be adjacent nodes in the Dynkin diagram.
Then the subsequence of $(i_1, \cdots, i_r)$ consisting of $i$ and $j$ 
is identical with the subsequence of $(j_1, \cdots, j_r)$ consisting of $i$ and $j$.
\item[(b)]
Let $i$ be a node in the Dynkin diagram.
Then the number of occurrence of $i$ in $(i_1, \cdots, i_r)$ 
is equal to the number of occurrence of $i$ in $(j_1, \cdots, j_r)$.
\end{enumerate}
\end{lemm}

We use these lemmas to prove the characterization of excited diagrams.

\begin{proof}[Proof of Proposition~\ref{prop:E=R}]
We follow the same idea used in the proof of \cite[Proposition~4.8]{GK}.
We denote by $\RR_P(F)$ the set of all subsets $D \subset P$ satisfying $\# D = \# F$ and $w_D = w_F$.

First we prove that $\EE_P(F) \subset \RR_P(F)$.
Since $F \in \RR_P(F)$, it is enough to show that, 
if $D' \in \EE_P(F)$ is obtained from $D \in \EE_P(F)$ by an elementary excitation, 
then $w_{D'} = w_D$.
Let $u$, $v \in P$ be elements such that $[v,u]$ is a $d_k$-interval, 
$[v,u] \cap D \cap N_{c(u)} = \emptyset$ and 
$D' = \alpha_u(D) = D \setminus \{ u \} \cup \{ v \}$.
If $v = p_k$, $u = p_l$ and $\{ j : k < j < l, \ p_j \in D \} = \{ j_1, \dots, j_m \}$ 
($j_1 < \dots < j_m$), then we have
$$
w_D = \cdots s(p_{j_1}) \cdots s(p_{j_m}) s(p_l) \cdots,
\quad
w_{D'} = \cdots s(p_k) s(p_{j_1}) \cdots s(p_{j_m}) \cdots.
$$
Then by using Lemma~\ref{lem:comm}, we have $w_{D'} = w_D$.

Next we prove that $\RR_P(F) \subset \EE_P(F)$.
Since $F$ is an order filter, we can take a linear extension of $P$ 
such that $F = \{ p_{n+1}, \cdots, p_N \}$, where $n = \# (P \setminus F)$.
We define the energy $e(D)$ of any subset $D \subset P$ with $\# D = \# F$ by putting
$$
e(D) = \sum_{p_k \in F} k - \sum_{p_k \in D} k.
$$
Then, by the assumption on our linear extension, 
we see that $e(D) \ge 0$ and that $e(D) = 0$ if and only if $D = F$.

We proceed by induction on $e(D)$ to prove $D \in \RR_P(F)$ implies $D \in \EE_P(F)$.
Let $D \in \RR_P(F)$.
If $e(D) = 0$, then we have $D = F \in \EE_P(F)$.

Suppose that $e(D) > 0$.
Then there exists an element $u \in F$ such that $u \not\in D$.
Let $u$ be the last element, \ie $u = p_l$ with $l$ largest, 
satisfying $u \in F$ and $u \not\in D$.
Since $w_D = w_F$ and $w_F$ is fully commutative, it follows from Lemma~\ref{lem:fc} (b)
that there exists an element $v \in D$ with the same color $i$ as $u$.
Let $v = p_k$ ($1 \le k < l$) be the last element satisfying $v \in D$ and $c(v) = c(u) = i$.
Then, by the choice of $k$ and $l$, we see that the reduced expressions of $w_F$ and $w_D$
are of the form $w_F = \cdots s_i s_{h_1} \cdots s_{h_r}$ ($s_i$ corresponds to $u$) 
and $w_D = \cdots s_i s_{j_1} \cdots s_{j_m} s_{h_1} \cdots s_{h_r}$ ($s_i$ corresponds to $v$) 
with no $i$ appearing in the segment $(j_1, \dots, j_m)$.
If $z \in [v,u] \cap D \cap N_i$, then $j = c(z)$ appears in the segment $(j_1, \dots, j_m)$ 
and this contradicts to the assertion of Lemma~\ref{lem:fc} (a).
Hence we have $[v,u] \cap D \cap N_i = \emptyset$.
If we put $D' = D \setminus \{ v \} \cup \{ u \}$,
then $w_{D'} = w_F$ by Lemma~\ref{lem:comm} and $e(D') < e(D)$.
Then by the induction hypothesis we have $D' \in \EE_P(F)$.
Since $D$ is obtained from $D'$ by a sequence of elementary excitations, 
we obtain $D \in \EE_P(F)$.
\end{proof}

Next we give a non-recursive description of the set of excited peaks $B(D)$, 
which implies that $B(D)$ is well-defined, \ie 
it is independent of the choices of elementary excitations to reach $D$ from $F$.

\begin{prop}
\label{prop:BD}
Let $D \in \EE_P(F)$ be an excited diagram.
\begin{enumerate}
\item[(a)]
The following are equivalent for $x \in P$:
\begin{enumerate}
\item[(i)]
$x \in B(D)$.
\item[(ii)]
There exists an element $y \in D_{c(x)}$ such that $y<x$ and $[y,x] \cap D \cap N_{c(x)} = \emptyset$.
\end{enumerate}
\item[(b)]
We have $D \cap B(D) = \emptyset$.
\end{enumerate}
\end{prop}

In the proof of this proposition, we utilize the following lemma, 
which will be used also in the sequel of this section.

\begin{lemm}
\label{lem:interval}
Let $x$, $y$, $z$, $u$ and $v$ be elements of $P$ such that $c(x) = c(y) = c(z) = i$, 
$c(u) = c(v) = j$,
$z<y<x$ and $v<u$.
Let $D$ be a subset of $P$.
Then we have
\begin{enumerate}
\item[(a)]
If $[z,y] \cap D \cap N_i = \emptyset$ and $[y,x] \cap D \cap N_i = \emptyset$, 
then we have $[z,x] \cap D \cap N_i = \emptyset$.
\item[(b)]
Suppose $u \in D$ and $v \not\in D$.
If $[y,x] \cap D \cap N_i = \emptyset$ and $x \not\in [v,u] \cap N_j$, 
then we have $[y,x] \cap (D \cup \{ v \}) \cap N_i = \emptyset$.
\item[(c)]
Suppose $u \in D$ and $v \not\in D$.
If $[y,x] \cap (D \setminus \{ u \} \cup \{ v \}) \cap N_i = \emptyset$ 
and $y \not\in [v,u] \cap N_j$, 
then we have $[y,x] \cap D \cap N_i = \emptyset$.
\item[(d)]
Suppose that $u \not\in D$ and $v \in D$.
If $[y,x] \cap D \cap N_i = \emptyset$ and $y \not\in [v,u] \cap N_j$, 
then we have $[y,x] \cap (D \cup \{ u \}) \cap N_i = \emptyset$.
\end{enumerate}
\end{lemm}

\begin{proof}
(a)
By using Property (C5) in Proposition~\ref{prop:dc_color}, 
we have $[z,x] \cap N_i = ( [z,y] \cap N_i ) \cup ( [y,x] \cap N_i)$.

(b)
Assume to the contrary that $[y,x] \cap (D \cup \{ v \}) \cap N_i \neq \emptyset$.
Since $[y,x] \cap D \cap N_i = \emptyset$, we have $v \in [y,x] \cap N_i$.
Thus $j = c(u) = c(v)$ is adjacent to $i = c(x)$ in $\Gamma$.
Since $u \in D$ and $[y,x] \cap D \cap N_i = \emptyset$, we have $u \not\in [y,x] \cap N_i$.
Hence, by using Property (C5) in Proposition~\ref{prop:dc_color},
we see that $y < v < x < u$ and $x \in [v,u] \cap N_j$, 
which contradicts to the assumption.

(c)
By an argument similar to (b), we can show that, if $[y,x] \cap D \cap N_i \neq \emptyset$, 
then $v < y < u < x$, which contradicts to $y \not\in [v,u] \cap N_j$.

(d)
By an argument similar to (b), we can show that, if $[y,x] \cap (D \cup \{ u \}) \cap N_i \neq \emptyset$, 
then $v < y < u < x$, which contradicts to $y \not\in [v,u] \cap N_j$.
\end{proof}

\begin{proof}[Proof of Proposition~\ref{prop:BD}]
We denote by $B'(D)$ the subset of $P$ consisting of elements $x \in P$ satisfying the condition (ii) in (a),
and prove $B(D) = B'(D)$ and $D \cap B'(D) = \emptyset$.
We proceed by induction on the number of elementary excitations to reach $D$ from $F$.

We begin with considering the case where $D = F$.
Let $x$ and $y$ be elements of $F$ with the same color $i$ satisfying $y<x$.
Then it follows from Properties (C4) and (C5) in Proposition~\ref{prop:dc_color} that 
an element $z$ covered by $x$ or covers $y$ belongs to $[y,x] \cap F \cap N_i$.
Hence we have $B'(F) = \emptyset = B(F)$ and $F \cap B'(F) = \emptyset$.

We prove that $B(\alpha_u(D)) = B'(\alpha_u(D))$ and $\alpha_u(D) \cap B'(\alpha_u(D)) = \emptyset$ 
for $D \in \EE_P(F)$ and a $D$-active element $u \in D$.
Let $v$ be the element such that $v<u$,
$c(v)=c(u)$, $[v,u] \cap D \cap N_{c(u)} = \emptyset$,
$[v,u]$ is a $d_k$-interval
and $\alpha_u(D) = D \setminus \{ u \} \cup \{ v \}$.
Recall that, by definition, 
$B(\alpha_u(D)) = \left( B(D) \setminus ([v,u] \cap N_{c(u)}) \right) \cup \{ u \}$.

First we show $B(\alpha_u(D)) \subset B'(\alpha_u(D))$.
Since $v \in \alpha_u(D)$ and $[v,u] \cap D \cap N_{c(u)} = \emptyset$, 
we have $u \in B'(\alpha_u(D))$.
Let $x \in B(D)$ such that $x \not\in [v,u]  \cap N_{c(u)}$.
Since $B(D) = B'(D)$ by the induction hypothesis, 
there exists $y \in D_{c(x)}$ such that $y<x$ and $[y,x] \cap D \cap N_{c(x)} = \emptyset$.
Then, by using Lemma~\ref{lem:interval} (b), 
we have $[y,x] \cap \alpha_u(D) \cap N_{c(x)} = \emptyset$, hence $x \in B'(\alpha_u(D))$.

Next, in order to show $B'(\alpha_u(D)) \subset B(\alpha_u(D))$, 
we take an element $x \in B'(\alpha_u(D))$ such that $x \neq u$ 
and prove $x \in B(D) \setminus ([v,u] \cap N_{c(u)})$.
Then there exists $y \in (\alpha_u(D))_{c(x)}$ such that 
$y<x$ and $[y,x] \cap \alpha_u(D) \cap N_{c(x)} = \emptyset$.
If $y = v$, then $\alpha_u(D)\cap N_{c(x)}=D\cap N_{c(x)}$, 
hence $[u,x] \cap D \cap N_{c(x)} \subset [v,x] \cap \alpha_u(D) \cap N_{c(x)} = \emptyset$.
Since $u \in D$, we have $x \in B'(D)=B(D)$ and $x \not\in N_{c(u)}$.
We consider the case where $y \neq v$.
In this case, $y \in D, y\not\in [v,u]\cap N_{c(u)}$ 
and it follows from Lemma~\ref{lem:interval} (c) 
that $[y,x] \cap D \cap N_{c(x)} = \emptyset$,
thus $x\in B'(D)=B(D)$.
Also we have $x \not\in [v,u] \cap N_{c(u)}$.
In fact, if $x \in [v,u] \cap N_{c(u)}$, 
then $c(y) = c(x)$ is adjacent to $c(u)$ 
and it follows from $y \in D$ and $[v,u] \cap D \cap N_{c(u)} = \emptyset$ that $y \not\in [v,u]\cap N_{c(u)}$.
Hence, by using Property (C5) in Proposition~\ref{prop:dc_color}, 
we have $y < v < x < u$ and 
this contradicts to $[y,x] \cap \alpha_u(D) \cap N_{c(x)} = \emptyset$.
Therefore we have $x \in B(D) \setminus ([v,u] \cap N_{c(u)})$.

Finally we show that $\alpha_u(D) \cap B'(\alpha_u(D)) = \emptyset$.
Let $x$ and $y$ be elements of $\alpha_u(D)$ with the same color $i$ satisfying $y<x$.
Since $\alpha_u(D) = D \setminus \{ u \} \cup \{ v\}$, it is enough to 
show 
$[y,x] \cap \alpha_u(D) \cap N_i \neq \emptyset$
for the following three cases:
$$
\text{Case~1. $x$, $y \in D$,}
\quad
\text{Case~2. $y = v$,}
\quad
\text{Case~3. $x = v$.}
$$
In Case~1, assume that $[y,x] \cap \alpha_u(D) \cap N_i = \emptyset$.
If $[y,x] \cap D \cap N_i = \emptyset$, then we have $x \in B(D)$ 
and this contradicts to the induction hypothesis $D \cap B(D) = \emptyset$.
Hence $[y,x] \cap D \cap N_i \neq \emptyset$ 
and this implies $u \in [y,x] \cap N_i$, 
\ie $y < u < x$ and $c(u) = c(v)$ is adjacent to $i$ in $\Gamma$.
Since $v \in \alpha_u(D)$, we have $v \not\in [y,x] \cap N_i$.
Therefore, by using Property (C5) in Proposition~\ref{prop:dc_color}, we see that 
$v < y < u < x$, which contradicts to the $D$-activity of $u$.
In Case~2, we have $v < u < x$ because $[v,u]$ is a $d_k$-interval 
and there is no element $w \in [v,u]$ with $w \neq v, u$ and $c(w) = c(v) = c(u)$ 
by Property (C2) in Proposition~\ref{prop:dc_color}.
Since $[u,x] \cap D \cap N_i \neq \emptyset$ by the induction hypothesis, 
we have $[v,x] \cap \alpha_u(D) \cap N_i \neq \emptyset$.
In Case~3, we have $[v,u] \cap D \cap N_i = \emptyset$ ($u$ is $D$-active) 
and $[y,u] \cap D \cap N_i \neq \emptyset$ (the induction hypothesis).
Hence by using Lemma~\ref{lem:interval} (a) we obtain $[y,v] \cap \alpha_u(D) \cap N_i \neq \emptyset$.
Therefore we see that any element satisfying the condition (ii) in (a) for $\alpha_u(D)$ 
does not belong to $\alpha_u(D)$, and $\alpha_u(D) \cap B'(\alpha_u(D)) = \emptyset$.
This completes the proof.
\end{proof}

\subsection{%
$K$-theoretical excited diagrams
}

We define $K$-theoretical excited diagrams and study their properties.
For shapes and shifted shapes, these diagrams were introduced in \cite{GK}.

\begin{defi}
\label{def:K-excited}
Let $P$ be a $d$-complete poset and let $F$ be an order filter of $P$.
\begin{enumerate}
\item[(a)]
Let $D$ be a subset of $P$ and $u \in D$.
If $u$ is $D$-active and $[v,u]$ is a $d_k$-interval, 
then we define $\alpha^*_u(D)$ to be the subset of $P$ obtained by 
adding $v$ to $D$.
We call this operation a \emph{$K$-theoretical elementary excitation}.
\item[(b)]
A \emph{$K$-theoretical excited diagram} of $F$ in $P$ is a subset of $P$ 
obtained from $F$ after a sequence of ordinary and $K$-theoretical elementary excitations on active elements. 
Let $\EE^*_P(F)$ be the set of all $K$-theoretical excited diagrams of $F$ in $P$.
\end{enumerate}
\end{defi}

\begin{exam}
\label{ex:shifted_K-excited}
If $P = S(5,4,2,1)$ is the shifted shape corresponding to a strict partition $(5,4,2,1)$ 
and $F = S(3,1)$, then there are $11$ $K$-theoretical excited diagrams in $\EE^{*}_P(F)$ 
shown in Figure~\ref{fig:shifted_K-excited}, 
and five of them are ordinary excited diagrams in $\EE_P(F)$.
In Figure~\ref{fig:shifted_K-excited}, 
the shaded cells form a $K$-theoretical exited diagram, 
and the arrow $D \longrightarrow D'$ (\resp $D \Longrightarrow D'$) indicates that 
$D'$ is obtained from $D$ by an ordinary (\resp $K$-theoretical) elementary excitation.
\begin{figure}[!htb]
$$
\setlength{\unitlength}{1.1pt}
\xymatrix@!C=9mm{
&&&
\raisebox{-20pt}{
\begin{picture}(50,40)
\fboxsep=0mm
\put(0,30){\colorbox[gray]{0.7}{\makebox(10,10){}}}
\put(10,30){\colorbox[gray]{0.7}{\makebox(10,10){}}}
\put(20,30){\colorbox[gray]{0.7}{\makebox(10,10){}}}
\put(10,20){\colorbox[gray]{0.7}{\makebox(10,10){}}}
\put(0,40){\line(1,0){50}}
\put(0,30){\line(1,0){50}}
\put(10,20){\line(1,0){40}}
\put(20,10){\line(1,0){20}}
\put(30,0){\line(1,0){10}}
\put(0,30){\line(0,1){10}}
\put(10,20){\line(0,1){20}}
\put(20,10){\line(0,1){30}}
\put(30,0){\line(0,1){40}}
\put(40,0){\line(0,1){40}}
\put(50,20){\line(0,1){20}}
\end{picture}
}
\ar[llld]
\ar[ld]
\ar@{=>}[rd]
\ar@{=>}[rrrd]
\\
\raisebox{-20pt}{
\begin{picture}(50,40)
\fboxsep=0mm
\put(0,30){\colorbox[gray]{0.7}{\makebox(10,10){}}}
\put(10,30){\colorbox[gray]{0.7}{\makebox(10,10){}}}
\put(30,20){\colorbox[gray]{0.7}{\makebox(10,10){}}}
\put(10,20){\colorbox[gray]{0.7}{\makebox(10,10){}}}
\put(0,40){\line(1,0){50}}
\put(0,30){\line(1,0){50}}
\put(10,20){\line(1,0){40}}
\put(20,10){\line(1,0){20}}
\put(30,0){\line(1,0){10}}
\put(0,30){\line(0,1){10}}
\put(10,20){\line(0,1){20}}
\put(20,10){\line(0,1){30}}
\put(30,0){\line(0,1){40}}
\put(40,0){\line(0,1){40}}
\put(50,20){\line(0,1){20}}
\end{picture}
}
\ar[d]
\ar@{=>}[rrd]
&&
\raisebox{-20pt}{
\begin{picture}(50,40)
\fboxsep=0mm
\put(0,30){\colorbox[gray]{0.7}{\makebox(10,10){}}}
\put(10,30){\colorbox[gray]{0.7}{\makebox(10,10){}}}
\put(20,30){\colorbox[gray]{0.7}{\makebox(10,10){}}}
\put(30,0){\colorbox[gray]{0.7}{\makebox(10,10){}}}
\put(0,40){\line(1,0){50}}
\put(0,30){\line(1,0){50}}
\put(10,20){\line(1,0){40}}
\put(20,10){\line(1,0){20}}
\put(30,0){\line(1,0){10}}
\put(0,30){\line(0,1){10}}
\put(10,20){\line(0,1){20}}
\put(20,10){\line(0,1){30}}
\put(30,0){\line(0,1){40}}
\put(40,0){\line(0,1){40}}
\put(50,20){\line(0,1){20}}
\end{picture}
}
\ar[lld]
\ar@{=>}[rrd]
&&
\raisebox{-20pt}{
\begin{picture}(50,40)
\fboxsep=0mm
\put(0,30){\colorbox[gray]{0.7}{\makebox(10,10){}}}
\put(10,30){\colorbox[gray]{0.7}{\makebox(10,10){}}}
\put(20,30){\colorbox[gray]{0.7}{\makebox(10,10){}}}
\put(30,20){\colorbox[gray]{0.7}{\makebox(10,10){}}}
\put(10,20){\colorbox[gray]{0.7}{\makebox(10,10){}}}
\put(0,40){\line(1,0){50}}
\put(0,30){\line(1,0){50}}
\put(10,20){\line(1,0){40}}
\put(20,10){\line(1,0){20}}
\put(30,0){\line(1,0){10}}
\put(0,30){\line(0,1){10}}
\put(10,20){\line(0,1){20}}
\put(20,10){\line(0,1){30}}
\put(30,0){\line(0,1){40}}
\put(40,0){\line(0,1){40}}
\put(50,20){\line(0,1){20}}
\end{picture}
}
\ar[d]
\ar@{=>}[rrd]
&&
\raisebox{-20pt}{
\begin{picture}(50,40)
\fboxsep=0mm
\put(0,30){\colorbox[gray]{0.7}{\makebox(10,10){}}}
\put(10,30){\colorbox[gray]{0.7}{\makebox(10,10){}}}
\put(20,30){\colorbox[gray]{0.7}{\makebox(10,10){}}}
\put(10,20){\colorbox[gray]{0.7}{\makebox(10,10){}}}
\put(30,00){\colorbox[gray]{0.7}{\makebox(10,10){}}}
\put(0,40){\line(1,0){50}}
\put(0,30){\line(1,0){50}}
\put(10,20){\line(1,0){40}}
\put(20,10){\line(1,0){20}}
\put(30,0){\line(1,0){10}}
\put(0,30){\line(0,1){10}}
\put(10,20){\line(0,1){20}}
\put(20,10){\line(0,1){30}}
\put(30,0){\line(0,1){40}}
\put(40,0){\line(0,1){40}}
\put(50,20){\line(0,1){20}}
\end{picture}
}
\ar[lllld]
\ar@{=>}[d]
\\
\raisebox{-20pt}{
\begin{picture}(50,40)
\fboxsep=0mm
\put(0,30){\colorbox[gray]{0.7}{\makebox(10,10){}}}
\put(10,30){\colorbox[gray]{0.7}{\makebox(10,10){}}}
\put(30,20){\colorbox[gray]{0.7}{\makebox(10,10){}}}
\put(30,00){\colorbox[gray]{0.7}{\makebox(10,10){}}}
\put(0,40){\line(1,0){50}}
\put(0,30){\line(1,0){50}}
\put(10,20){\line(1,0){40}}
\put(20,10){\line(1,0){20}}
\put(30,0){\line(1,0){10}}
\put(0,30){\line(0,1){10}}
\put(10,20){\line(0,1){20}}
\put(20,10){\line(0,1){30}}
\put(30,0){\line(0,1){40}}
\put(40,0){\line(0,1){40}}
\put(50,20){\line(0,1){20}}
\end{picture}
}
\ar[d]
\ar@{=>}[rrd]
&&
\raisebox{-20pt}{
\begin{picture}(50,40)
\fboxsep=0mm
\put(0,30){\colorbox[gray]{0.7}{\makebox(10,10){}}}
\put(10,30){\colorbox[gray]{0.7}{\makebox(10,10){}}}
\put(10,20){\colorbox[gray]{0.7}{\makebox(10,10){}}}
\put(30,20){\colorbox[gray]{0.7}{\makebox(10,10){}}}
\put(30,0){\colorbox[gray]{0.7}{\makebox(10,10){}}}
\put(0,40){\line(1,0){50}}
\put(0,30){\line(1,0){50}}
\put(10,20){\line(1,0){40}}
\put(20,10){\line(1,0){20}}
\put(30,0){\line(1,0){10}}
\put(0,30){\line(0,1){10}}
\put(10,20){\line(0,1){20}}
\put(20,10){\line(0,1){30}}
\put(30,0){\line(0,1){40}}
\put(40,0){\line(0,1){40}}
\put(50,20){\line(0,1){20}}
\end{picture}
}
&&
\raisebox{-20pt}{
\begin{picture}(50,40)
\fboxsep=0mm
\put(0,30){\colorbox[gray]{0.7}{\makebox(10,10){}}}
\put(10,30){\colorbox[gray]{0.7}{\makebox(10,10){}}}
\put(20,30){\colorbox[gray]{0.7}{\makebox(10,10){}}}
\put(30,20){\colorbox[gray]{0.7}{\makebox(10,10){}}}
\put(30,0){\colorbox[gray]{0.7}{\makebox(10,10){}}}
\put(0,40){\line(1,0){50}}
\put(0,30){\line(1,0){50}}
\put(10,20){\line(1,0){40}}
\put(20,10){\line(1,0){20}}
\put(30,0){\line(1,0){10}}
\put(0,30){\line(0,1){10}}
\put(10,20){\line(0,1){20}}
\put(20,10){\line(0,1){30}}
\put(30,0){\line(0,1){40}}
\put(40,0){\line(0,1){40}}
\put(50,20){\line(0,1){20}}
\end{picture}
}
&&
\raisebox{-20pt}{
\begin{picture}(50,40)
\fboxsep=0mm
\put(0,30){\colorbox[gray]{0.7}{\makebox(10,10){}}}
\put(10,30){\colorbox[gray]{0.7}{\makebox(10,10){}}}
\put(30,20){\colorbox[gray]{0.7}{\makebox(10,10){}}}
\put(10,20){\colorbox[gray]{0.7}{\makebox(10,10){}}}
\put(20,30){\colorbox[gray]{0.7}{\makebox(10,10){}}}
\put(30,0){\colorbox[gray]{0.7}{\makebox(10,10){}}}
\put(0,40){\line(1,0){50}}
\put(0,30){\line(1,0){50}}
\put(10,20){\line(1,0){40}}
\put(20,10){\line(1,0){20}}
\put(30,0){\line(1,0){10}}
\put(0,30){\line(0,1){10}}
\put(10,20){\line(0,1){20}}
\put(20,10){\line(0,1){30}}
\put(30,0){\line(0,1){40}}
\put(40,0){\line(0,1){40}}
\put(50,20){\line(0,1){20}}
\end{picture}
}
\\
\raisebox{-20pt}{
\begin{picture}(50,40)
\fboxsep=0mm
\put(0,30){\colorbox[gray]{0.7}{\makebox(10,10){}}}
\put(20,20){\colorbox[gray]{0.7}{\makebox(10,10){}}}
\put(30,20){\colorbox[gray]{0.7}{\makebox(10,10){}}}
\put(30,0){\colorbox[gray]{0.7}{\makebox(10,10){}}}
\put(0,40){\line(1,0){50}}
\put(0,30){\line(1,0){50}}
\put(10,20){\line(1,0){40}}
\put(20,10){\line(1,0){20}}
\put(30,0){\line(1,0){10}}
\put(0,30){\line(0,1){10}}
\put(10,20){\line(0,1){20}}
\put(20,10){\line(0,1){30}}
\put(30,0){\line(0,1){40}}
\put(40,0){\line(0,1){40}}
\put(50,20){\line(0,1){20}}
\end{picture}
}
&&
\raisebox{-20pt}{
\begin{picture}(50,40)
\fboxsep=0mm
\put(0,30){\colorbox[gray]{0.7}{\makebox(10,10){}}}
\put(10,30){\colorbox[gray]{0.7}{\makebox(10,10){}}}
\put(20,20){\colorbox[gray]{0.7}{\makebox(10,10){}}}
\put(30,20){\colorbox[gray]{0.7}{\makebox(10,10){}}}
\put(30,0){\colorbox[gray]{0.7}{\makebox(10,10){}}}
\put(0,40){\line(1,0){50}}
\put(0,30){\line(1,0){50}}
\put(10,20){\line(1,0){40}}
\put(20,10){\line(1,0){20}}
\put(30,0){\line(1,0){10}}
\put(0,30){\line(0,1){10}}
\put(10,20){\line(0,1){20}}
\put(20,10){\line(0,1){30}}
\put(30,0){\line(0,1){40}}
\put(40,0){\line(0,1){40}}
\put(50,20){\line(0,1){20}}
\end{picture}
}
}
$$
\caption{$K$-theoretical excited diagrams of $S(3,1)$ in $S(5,4,2,1)$}
\label{fig:shifted_K-excited}
\end{figure}
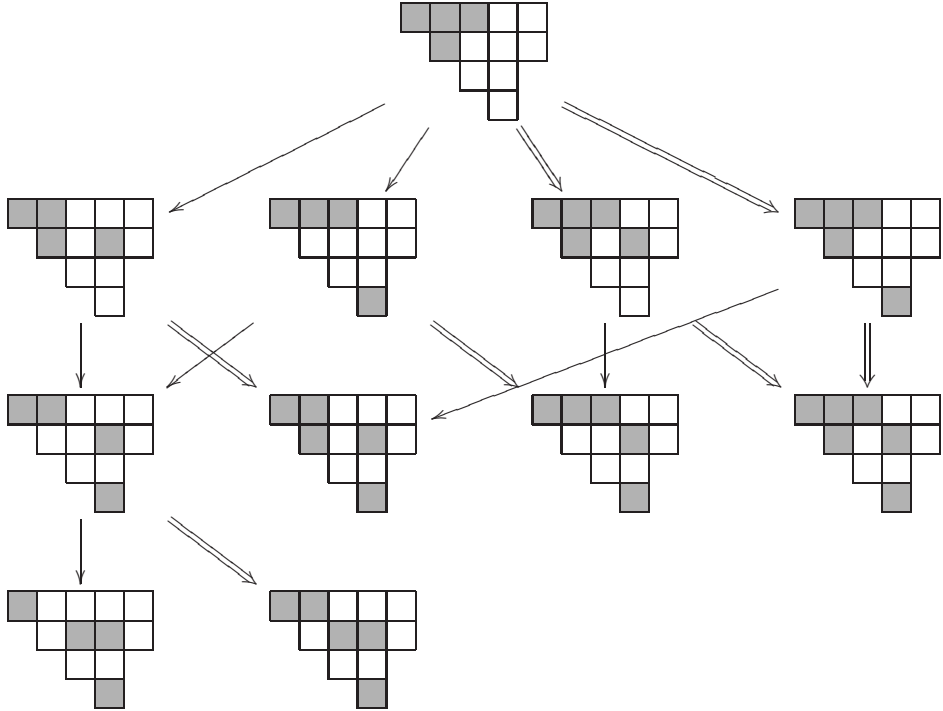
\end{exam}

For a fixed linear extension of $P$ and 
a subset $D = \{ p_{i_1}, \cdots, p_{i_r} \}$ ($i_1 < \cdots < i_r$) of $P$, 
we define an element $w^*_D \in W$ by putting
$$
w^*_D = s(p_{i_1}) * s(p_{i_2}) * \cdots * s(p_{i_r}),
$$
where $* : W \times W \to W$ is the associative product, called the Demazure product, defined by
$$
s_i * w
 =
\begin{cases}
 s_i w &\text{if $l(s_i w) = l(w) + 1$,} \\
 w &\text{if $l(s_i w) = l(w) - 1$.}
\end{cases}
$$
For the fundamental properties of the Demazure product
(also called the Hecke product), we refer to \cite[Section 3]{BM}.
Since $w_P$ is fully commutative (Proposition~\ref{prop:w_P}), the element $w^*_D$ is 
independent of the choices of linear extensions of $P$.

The following proposition is a key to rephrase the Billey-type formula for equivariant $K$-theory 
in terms of combinatorics of $d$-complete posets (see Proposition~\ref{prop:Billey_P}).

\begin{prop}
\label{prop:E*=R*}
Let $P$ be a connected $d$-complete poset and $F$ an order filter of $P$.
Then a subset $D \subset P$ is a $K$-theoretical excited diagram of $F$ in $P$ 
if and only if $w^*_D = w_F$.
\end{prop}

We follow the same idea as the proof of \cite[Proposition~4.8]{GK}.

\begin{lemm}
\label{lem:red}
\textup{(\cite[Lemma~3.1 and Proposition~3.4]{GK})}
Let $v \in W$ and $v = s_{i_1} * \cdots * s_{i_r}$ with $i_1, \dots, i_r \in I$.
\begin{enumerate}
\item[(a)]
There is a increasing sequence $1 \le k_1 < k_2 < \dots < k_l \le r$ such that 
$v = s_{i_{k_1}} s_{i_{k_2}} \dots s_{i_{k_l}}$ is a reduced expression of $v$.
In particular, we have $l(v) \le r$.
\item[(b)]
If $l(v) = r$, then $v = s_{i_1} \cdots s_{i_r}$.
\item[(c)]
If $v$ is fully commutative and $l(v) < r$, 
then there exist $a < b$ such that $s_{i_a} = s_{i_b}$ and
$s_{i_a}$ commutes with $s_{i_c}$ for every $a < c < b$.
\end{enumerate}
\end{lemm}

\begin{proof}[Proof of Proposition~\ref{prop:E*=R*}]
We denote by $\RR^*_P(F)$ the set of all subsets $D \subset P$ satisfying $w^*_D = w_F$.

First we prove $\EE^*_P(F) \subset \RR^*_P(F)$.
Since $F \in \RR^*_P(F)$, it is enough to show that, 
if $D' \in \EE^*_P(F)$ is obtained from $D \in \EE^*_P(F)$ by an ordinary or $K$-theoretical elementary excitation, 
then $w^*_{D'} = w^*_D$.
Let $u$ be a $D$-active element and $[v,u]$ be the $d_k$-interval with top element $u$.
Let $v = p_k$, $u = p_l$ and $\{ j : k < j < l, p_j \in D \} = \{ j_1, \dots, j_m \}$ ($j_1 < \dots < j_m$).
If $D' = \alpha_u(D) = D \setminus \{ u \} \cup \{ v \}$, then we have
\begin{align*}
w^*_D &= \cdots * s(p_{j_1}) * \cdots * s(p_{j_m}) * s(p_l) * \cdots,
\\
w^*_{D'} &= \cdots * s(p_k) * s(p_{j_1}) * \cdots * s(p_{j_m}) * \cdots.
\end{align*}
If $D' = \alpha^*_u(D) = D \cup \{ v \}$, then we have
\begin{align*}
w^*_D &= \cdots * s(p_{j_1}) * \cdots * s(p_{j_m}) * s(p_l) * \cdots,
\\
w^*_{D'} &= \cdots * s(p_k) * s(p_{j_1}) * \cdots * s(p_{j_m}) * s(p_l) * \cdots.
\end{align*}
Since $u$ is $D$-active, it follows from Lemma~\ref{lem:comm} and $s(p_k)*s(p_l) = s(p_l)*s(p_l) = s(p_l)$ 
that $w^*_{D'} = w^*_D$ in both cases.

Next we prove $\RR^*_P(F) \subset \EE^*_P(F)$.
We proceed by induction on $\# D$ to prove $D \in \RR^*_P(F)$ implies $D \in \EE^*_P(F)$.
Since $w^*_D = w_F$, we have $\# D \ge l(w_F) = \# F$ by Lemma~\ref{lem:red} (a).
If $\# D = \# F$, then $w^*_D = w_D$ by Lemma~\ref{lem:red} (b), 
thus we have $D \in \EE_P(F) \subset \EE^*_P(F)$ by Proposition~\ref{prop:E=R}.
If $\# D > \# F$, then it follows from Lemma~\ref{lem:red} (c) that 
there exist $u$, $v \in D$ with $v < u$ such that $c(u) = c(v)$ and 
$s(u) = s(v)$ commutes with $s(w)$ for
every element $w$ between $u$ and $v$ in the expression of $w^*_D$.
If we put $D' = D \setminus \{ v \}$, then we see that $w^*_{D'} = w^*_D$.
Hence by the induction hypothesis we have $D' \in \EE^*_P(F)$.
Since $D$ is obtained from $D'$ by a sequence of ordinary and $K$-theoretical elementary excitations, 
we obtain $D \in \EE^*_P(F)$.
\end{proof}

The following proposition plays a crucial role in the proof of our main theorem 
(see the proof of Theorem~\ref{thm:xi2}).

\begin{prop}
\label{prop:E*}
Let $P$ be a connected $d$-complete poset and $F$ an order filter of $P$.
Then we have
\begin{equation}
\label{eq:E*}
\EE^*_P(F) = \bigsqcup_{D \in \EE_P(F)} \{ D \sqcup S : S \subset B(D) \}.
\end{equation}
\end{prop}

In order to prove this proposition, we prepare several lemmas.

\begin{lemm}
\label{lem:interval_BD}
Let $x$ and $y$ be elements of $P$ such that $y<x$ and $c(x) = c(y) = i$, 
and let $D\in \EE_P(F)$ be an excited diagram.
If $[y,x] \cap D \cap N_i = \emptyset$ and $y \in D$, 
then we have $[y,x] \cap B(D) \cap N_i = \emptyset$.
\end{lemm}

\begin{proof}
Assume to the contrary that $[y,x] \cap B(D) \cap N_i \neq \emptyset$ 
and take an element $z \in [y,x] \cap B(D) \cap N_i$.
By Proposition~\ref{prop:BD} (a), there exists $w \in D_{c(z)}$ such that 
$w < z$ and $[w,z] \cap D \cap N_{c(z)} = \emptyset$.
Then $c(w) = c(z)$ is adjacent to $i = c(y)$ in $\Gamma$ and $w \in D \cap N_{c(x)}$.
Since $[y,x] \cap D \cap N_{c(x)} = \emptyset$, we have $w \not\in [y,x]$.
Hence by using Property (C5) in Proposition~\ref{prop:dc_color} 
we see that $w < y < z < x$, which contracts to $[w,z] \cap D \cap N_{c(z)} = \emptyset$.
\end{proof}

For $E \in \EE^*_P(F)$, we define a subset $S(E)$ of $E$ by putting
\begin{equation}
\label{eq:def_SE}
S(E)
 = 
\left\{
x \in E :
\begin{matrix}
\text{there exists $y \in E_{c(x)}$ such that } \\
\text{$y<x$ and $[y,x] \cap E \cap N_{c(x)} = \emptyset$}
\end{matrix}
\right\}.
\end{equation}
It follows from Property (C4) in Proposition~\ref{prop:dc_color} that $S(F) = \emptyset$ 
for an order filter $F$ of $P$.

\begin{lemm}
\label{lem:SE1}
Let $E \in \EE^*_P(F)$ and $u \in E$ an $E$-active element.
Then we have
\begin{align}
\label{eq:SE11}
S(\alpha_u(E))
 &=
\begin{cases}
 S(E) \setminus \{ u \} \cup \{ v \} &\text{if $u \in S(E)$,} \\
 S(E) &\text{if $u \not\in S(E)$,}
\end{cases}
\\
\label{eq:SE12}
S(\alpha^*_u(E))
 &=
\begin{cases}
 S(E) \cup \{ v \} &\text{if $u \in S(E)$,} \\
 S(E) \cup \{ u \} &\text{if $u \not\in S(E)$,}
\end{cases}
\end{align}
where $[v,u]$ is the $d_k$-interval with top element $u$.
In particular, we have
$$
\# S(\alpha_u(E)) = \# S(E),
\quad
\# S(\alpha^*_u(E)) = \# S(E) + 1,
$$
and $S(E) = \emptyset$ if and only if $E \in \EE_P(F)$.
\end{lemm}

\begin{proof}
Since $\alpha_u(E) = E \setminus \{ u \} \cup \{ v \}$, the equality \eqref{eq:SE11} 
follows from the following three claims:
\begin{enumerate}
\item[(i)]
$S(\alpha_u(E)) \setminus \{ v \} \subset S(E)$.
\item[(ii)]
$u \in S(E)$ if and only if $v \in S(\alpha_u(E))$.
\item[(iii)]
$S(E) \subset S(\alpha_u(E))$.
\end{enumerate}

To prove (i), we take $x \in S(\alpha_u(E))$ such that $x \neq v$.
Then, by the definition \eqref{eq:def_SE}, there exists $y \in (\alpha_u(E))_{c(x)}$ such that $y<x$ 
and $[y,x] \cap \alpha_u(E) \cap N_{c(x)} = \emptyset$.
If $y = v$, then $[u,x] \cap E \cap N_{c(x)} \subset [v,x] \cap E \cap N_{c(x)} = \emptyset$, 
hence $x \in S(E)$.
We consider the case where $y \neq v$.
In this case, $y \in E$ and it follows from $[v,u] \cap E \cap N_{c(u)} = \emptyset$ 
that $y \not\in [v,u] \cap N_{c(u)}$.
Now we can use Lemma~\ref{lem:interval} (c) to obtain $[y,x] \cap E \cap N_{c(x)} = \emptyset$, 
hence $x \in S(E)$.

Next we prove (ii).
If $v\in S(\alpha_u(E))$, then there exists $w \in \alpha_u(E)_{c(v)}$ such that 
$w<v$ and $[w,v] \cap \alpha_u(E) \cap N_{c(v)} = \emptyset$.
Then by Lemma~\ref{lem:interval} (a) we have $[w,u] \cap E \cap N_{c(u)} = \emptyset$, 
hence $u \in S(E)$.
Conversely, 
if $u \in S(E)$, there exists $z \in E_{c(u)}$ such that $z < u$ and $[z,u] \cap E \cap N_{c(u)} = \emptyset$.
Then we have $[z,v] \cap \alpha_u(E) \cap N_{c(v)} \subset [z,u] \cap \alpha_u(E) \cap N_{c(u)} = \emptyset$, 
hence $v \in S(\alpha_u(E))$.

To prove (iii), we take $x \in S(E)$.
By the definition \eqref{eq:def_SE}, 
there exists $y \in E_{c(x)}$ such that $y<x$ and $[y,x] \cap E \cap N_{c(x)} = \emptyset$.
Then we have $y = u$ or $y \in \alpha_u(E)$.
If $y=u$, then $[u,x] \cap E \cap N_{c(x)} = \emptyset$ and 
$[v,u] \cap \alpha_u(E) \cap N_{c(u)} = \emptyset$ by the $E$-activity of $u$.
Hence by using Lemma~\ref{lem:interval} (a) we have $[v,x] \cap \alpha_u(E) \cap N_{c(x)} = \emptyset$ 
and $x \in S(\alpha_u(E))$.
We consider the case where $y \in \alpha_u(E)$.
Since $[v,u] \cap E \cap N_{c(u)} = \emptyset$ and $x \in E$, we have $x \not\in [v,u] \cap N_{c(u)}$.
Hence by using Lemma~\ref{lem:interval} (b) we see $[y,x] \cap \alpha_u(E) \cap N_{c(x)} = \emptyset$ 
and $x \in S(\alpha_u(E))$.

Since $\alpha^*_u(E) = E \cup \{ v \}$ and 
$u\in S(\alpha^*_u(E))$, the equality \eqref{eq:SE12} 
follows from the following three claims:
\begin{enumerate}
\item[(iv)]
$S(\alpha^*_u(E)) \setminus \{ u, v \} \subset S(E)$.
\item[(v)]
$u \in S(E)$ if and only if $v \in S(\alpha^*_u(E))$.
\item[(vi)]
$S(E) \subset S(\alpha^*_u(E))$.
\end{enumerate}

To prove (iv), we take $x \in S(\alpha^*_u(E))$ such that $x \neq u$, $v$.
Then there exists $y \in (\alpha^*_u(E))_{c(x)}$ such that $y<x$ 
and $[y,x] \cap \alpha^*_u(E) \cap N_{c(x)} = \emptyset$.
If $y = v$, then $[u,x] \cap E \cap N_{c(x)} \subset [v,x] \cap E \cap N_{c(x)} = \emptyset$, 
hence $x \in S(E)$.
If $y \in E$, then $[y,x] \cap E \cap N_{c(x)} \subset [y,x] \cap \alpha^*_u(E) \cap N_{c(x)} = \emptyset$, 
hence $x \in S(E)$.

Next we prove (v).
If $v \in S(\alpha^{*}_u(E))$, then there exists $w \in \alpha^*_u(E)_{c(v)}$ such that
$w<v$ and $[w,v] \cap \alpha^*_u(E) \cap N_{c(v)} = \emptyset$.
Then by using Lemma~\ref{lem:interval} (a) we have $[w,u] \cap E \cap N_{c(u)} = \emptyset$, 
hence $u \in S(E)$.
Conversely, if $u \in S(E)$, there exists $z \in E_{c(u)}$ such that $z<u$ and 
$[z,u] \cap E \cap N_{c(u)} = \emptyset$.
Then we have $[z,v] \cap \alpha^*_u(E) \cap N_{c(v)} \subset [z,u] \cap \alpha^*_u(E) \cap N_{c(v)} = \emptyset$, 
hence $v \in S(\alpha^*_u(E))$.

To prove (vi), we take $x \in S(E)$.
Then there exists $y \in E_{c(x)}$ such that $y < x$ and $[y,x] \cap E \cap N_{c(x)} = \emptyset$.
Since $[v,u] \cap E \cap N_{c(u)} = \emptyset$ and $x \in E$, we have $x \not\in [v,u] \cap N_{c(u)}$.
Hence by using Lemma~\ref{lem:interval} (b) we see $[y,x] \cap \alpha^*_u(E) \cap N_{c(x)} = \emptyset$ 
and $x \in S(\alpha^*_u(E))$.
\end{proof}

\begin{lemm}
\label{lem:SE2}
If $D \in \EE_P(F)$ and $S \subset B(D)$, 
then we have $D \sqcup S \in \EE^*_P(F)$ and $S = S(D \sqcup S)$.
In particular, if $E \in \EE^*_P(F)$ is expressed as $E = D \sqcup S = D' \sqcup S'$ 
with $D$, $D' \in \EE_P(F)$ and $S \subset B(D)$, $S' \subset B(D')$, then we have 
$D = D'$ and $S = S'$.
\end{lemm}

\begin{proof}
First we proceed by induction on $\# S$ to prove $D \sqcup S \in \EE^*_P(F)$.
If $S = \emptyset$, then we have $D \in \EE_P(F) \subset \EE^*_P(F)$.
If $S \neq \emptyset$, we take an element $x \in S$ and put $S' = S \setminus \{ x \}$.
By the induction hypothesis and Proposition~\ref{prop:E*=R*}, 
we have $D \sqcup S' \in \EE^*_P(F)$ and $w^*_{D \sqcup S'} = w_F$.
Using Proposition~\ref{prop:BD} (a), we see that there exists $y \in D_{c(x)}$ 
such that $y<x$ and $[y,x] \cap D \cap N_{c(x)} = \emptyset$.
Then it follows from Lemma~\ref{lem:interval_BD} that $[y,x] \cap (D \sqcup S) \cap N_{c(x)} = \emptyset$.
If $y = p_k$, $x = p_l$ and $\{ j : k < j < l, \ p_j \in D \sqcup S \} = \{ j_1, \dots, j_m \}$ 
($j_1 < \dots < j_m$), then we have 
\begin{align*}
w^*_{D \sqcup S} &= \cdots * s(p_k) * s(p_{j_1}) * \cdots * s(p_{j_m}) * s(p_l) * \cdots,
\\
w^*_{D \sqcup S'} &= \cdots * s(p_k) * s(p_{j_1}) * \cdots * s(p_{j_m}) * \cdots.
\end{align*}
By using Lemma~\ref{lem:comm} and $s(p_k)*s(p_l) = s(p_k) * s(p_k) = s(p_k)$, 
we obtain $w^*_{D \sqcup S} = w^*_{D \sqcup S'} = w_F$.
Hence by Proposition~\ref{prop:E*=R*} we have $D \sqcup S \in \EE^*_P(F)$.

Next we put $E = D \sqcup S$ and prove that $S = S(E)$.
In order to show the inclusion $S \subset S(E)$, we take $x \in S \subset B(D)$.
Then by Proposition~\ref{prop:BD} (a), there exists $y \in D_{c(x)}$ 
such that $[y,x] \cap D \cap N_{c(x)} = \emptyset$.
Hence by using Lemma~\ref{lem:interval_BD} we see that $[y,x] \cap E \cap N_{c(x)} = \emptyset$ 
and $x \in S(E)$.
In order to show the reverse inclusion $S(E) \subset S$, we take $x \in S(E)$ and prove $x \in B(D)$.
By the definition \eqref{eq:def_SE}, there exists $y \in E_{c(x)}$ 
such that $y < x$ and $[y,x] \cap E \cap N_{c(x)} = \emptyset$.
Since $D \subset E$, we have $[y,x] \cap D \cap N_{c(x)} = \emptyset$.
If $y \in D$, then by Proposition~\ref{prop:BD} (a) we have $x \in B(D)$.
If $y \in S$, then there exists $z \in D_{c(y)}$ such that $z<y$ and $[z,y] \cap D \cap N_{c(y)} = \emptyset$.
Then by using Lemma~\ref{lem:interval} (a) we have $[z,x] \cap D \cap N_{c(x)} = \emptyset$ and $x \in B(D)$.
\end{proof}

\begin{lemm}
\label{lem:SE3}
Let $E \in \EE^*_P(F)$ and $z \in S(E)$. If we put $E' = E \setminus \{ z \}$, then we have
\begin{enumerate}
\item[(a)]
$E' \in \EE^*_P(F)$.
\item[(b)]
$S( E' ) = S(E) \setminus \{ z \}$.
\end{enumerate}
\end{lemm}

\begin{proof}
By the definition \eqref{eq:def_SE}, 
there exists $w \in E_{c(z)}$ such that $w < z$ and $[w,z] \cap E \cap N_{c(z)} = \emptyset$.

(a)
If $w = p_k$, $z = p_l$ and $\{ j : k < j < l, \ p_j \in E \} = \{ j_1, \dots, j_m \}$ 
($j_1 < \dots < j_m$), then we have 
\begin{align*}
w^*_E &= \cdots * s(p_k) * s(p_{j_1}) * \cdots * s(p_{j_m}) * s(p_l) * \cdots,
\\
w^*_{E'} &= \cdots * s(p_k) * s(p_{j_1}) * \cdots * s(p_{j_m}) * \cdots.
\end{align*}
By using Lemma~\ref{lem:comm} and $s(p_k)*s(p_l) = s(p_k)$, 
we obtain $w^*_E = w^*_{E'}$.
Hence it follows from Proposition~\ref{prop:E*=R*} that $E' \in \EE^*_P(F)$.

(b)
First we prove that $S(E') \subset S(E)$.
Let $x \in S(E')$.
Then there exists $y \in (E')_{c(x)}$ such that $y < x $ and $[y,x] \cap E' \cap N_{c(x)} = \emptyset$.
Since $y \in E$ and $[w,z] \cap E \cap N_{c(z)} = \emptyset$, we have $y \not\in [w,z] \cap N_{c(z)}$.
Hence by using Lemma~\ref{lem:interval} (d) we have $[y,x] \cap E \cap N_{c(x)} = \emptyset$ and $x \in S(E)$.

Next we prove that $S(E) \setminus \{ z \} \subset S(E')$.
We take an element $x \in S(E)$ such that $x \neq z$.
Then there exists $y \in E_{c(x)}$ such that $y<x$ and $[y,x] \cap E \cap N_{c(x)} = \emptyset$.
If $y \neq z$, then $y \in E'$ and $[y,x] \cap E' \cap N_{c(x)} = \emptyset$, hence $x \in S(E')$.
We consider the case where $y = z$.
In this case $[z,x] \cap E' \cap N_{c(x)} = \emptyset$.
Since $[w,z] \cap E' \cap N_{c(z)} = \emptyset$, it follows from Lemma~\ref{lem:interval} (a) 
that $[w,x] \cap E' \cap N_{c(x)} = \emptyset$ and $x \in S(E')$.
\end{proof}

Now we are in position to give a proof of Proposition~\ref{prop:E*}.

\begin{proof}[Proof of Proposition~\ref{prop:E*}]
By using Lemma~\ref{lem:SE2}, it is enough to show that any $E \in \EE^*_P(F)$ 
can be written as $E = D \sqcup S$ with $D \in \EE_P(F)$ and $S \subset B(D)$.
Given $E \in \EE^*_P(F)$, we put $D = E \setminus S(E)$ and prove that 
$D \in \EE_P(F)$ and $S(E) \subset B(D)$.
We proceed by induction on $\# S(E)$.
If $S(E) = \emptyset$, then $E \in \EE_P(F)$ by Lemma~\ref{lem:SE1}.
We consider the case where $S(E) \neq \emptyset$.
Then we take $z \in S(E)$ and put $E' = E \setminus \{ z \}$.
Since $S(E') = S(E) \setminus \{ z \}$ by Lemma~\ref{lem:SE3} (b), 
we see that $D = E \setminus S(E) = E' \setminus S(E')$.
By the induction hypothesis, $D \in \EE_P(F)$ and $S(E') \subset B(D)$.
It remains to show that $z \in B(D)$.
By the definition \eqref{eq:def_SE}, 
there exists $w \in E_{c(z)}$ such that $w<z$ and $[w,z] \cap E \cap N_{c(z)} = \emptyset$.
Let $w$ be the minimal such element.
If $w \not\in D$, \ie $w \in S(E)$, 
then by definition there exists $w' \in E_{c(w)}$ such that $[w',w] \cap E \cap N_{c(w)} = \emptyset$.
Then it follows from Lemma~\ref{lem:interval} (a) that $[w',z] \cap E \cap N_{c(z)} = \emptyset$, 
which contradicts to the minimality of $w$.
Therefore we have $w \in D$ and $z \in B(D)$.
\end{proof}

\begin{rema}
\label{rem:heap2}
We can define the notion of $D$-active elements and ordinary and $K$-theoretical elementary excitations 
for a general dominant minuscule heap $H(w)$ just by replacing $d_k$-intervals with 
intervals $[v,u]$ such that $c(u) = c(v)$ and $[v,u] \cap H(w)_{c(u)} = \{ u, v \}$.
Then the arguments in this section work as well for $H(w)$.
In particular, Propositions~\ref{prop:E*=R*} and \ref{prop:E*} holds for $H(w)$.
\end{rema}

\begin{exam}
\label{ex:shited-B-excited}
Let $\mu$ be a strict partition and regard the shifted shape $S(\mu)$ as a heap for the Weyl group of type $B$ 
(see Example~\ref{ex:shifted-B}).
Then the notion of active elements, (ordinary and $K$-theoretical) elementary excitations and excited peaks 
are modified as follows.
Let $D$ be a subset of $S(\mu)$.
\begin{enumerate}
\item[(a)]
An element $u = (i,j) \in D$ is $D$-active if either
\begin{gather*}
\text{$i < j$ and $(i,j+1)$, $(i+1,j)$, $(i+1,j+1) \in S(\mu) \setminus D$, or}
\\
\text{$i = j$ and $(i,i+1)$, $(i+1,i+1) \in S(\mu) \setminus D$.}
\end{gather*}
\item[(b)]
If $u = (i,j) \in D$ is $D$-active, then we define an ordinary and $K$-theoretical 
elementary excitation by putting
\begin{align*}
\alpha_u(D) &= D \setminus \{ (i,j) \} \cup \{ (i+1,j+1) \},
\\
\alpha^*_u(D) &= D \cup \{ (i+1,j+1) \},
\end{align*}
respectively.
\item[(c)]
If $u = (i,j) \in D$ is $D$-active, then we define
$$
B(\alpha_u(D))
 =
\begin{cases}
 B(D) \setminus \{ (i,j+1), (i+1,j) \} \cup \{ (i,j) \} 
 &\text{if $i<j$,}
\\
 B(D) \setminus \{ (i,j+1) \} \cup \{ (i,j) \} 
 &\text{if $i=j$.}
\end{cases}
$$
\end{enumerate}
This notion of excited diagrams is the same as Ikeda--Naruse's excited diagrams of type $I$ 
introduced in \cite{IN1}.

For example, if $P = S(5,4,2,1)$ and $F = S(3,1)$, 
then there are $10$ excited diagrams of $F$ in $P$ as a heap for the type $B$ Weyl group.
See Figure~\ref{fig:shifted-B_excited}.
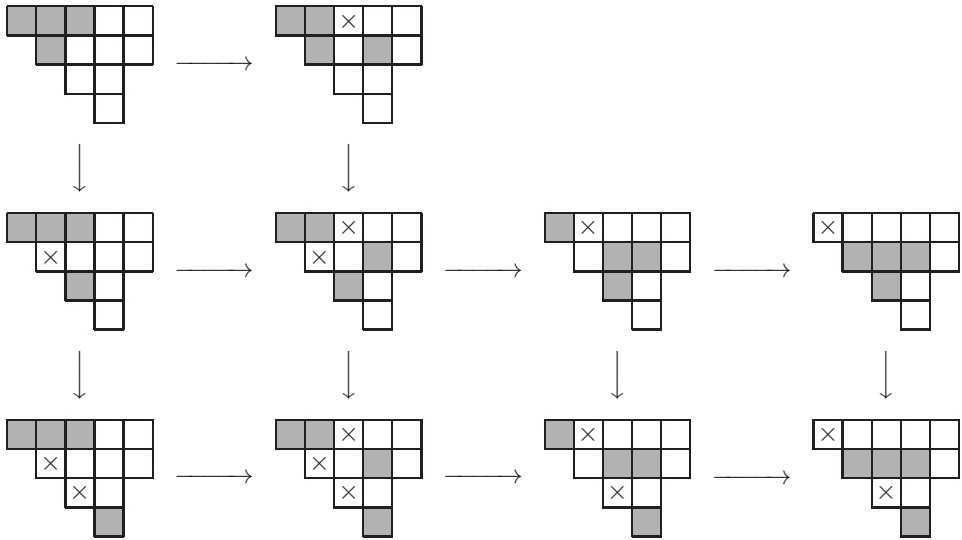
\begin{figure}[!htb]
$$
\setlength{\unitlength}{1.1pt}
\begin{CD}
\raisebox{-20pt}{
\begin{picture}(50,40)
\fboxsep=0mm
\put(0,30){\colorbox[gray]{0.7}{\makebox(10,10){}}}
\put(10,30){\colorbox[gray]{0.7}{\makebox(10,10){}}}
\put(20,30){\colorbox[gray]{0.7}{\makebox(10,10){}}}
\put(10,20){\colorbox[gray]{0.7}{\makebox(10,10){}}}
\put(0,40){\line(1,0){50}}
\put(0,30){\line(1,0){50}}
\put(10,20){\line(1,0){40}}
\put(20,10){\line(1,0){20}}
\put(30,0){\line(1,0){10}}
\put(0,30){\line(0,1){10}}
\put(10,20){\line(0,1){20}}
\put(20,10){\line(0,1){30}}
\put(30,0){\line(0,1){40}}
\put(40,0){\line(0,1){40}}
\put(50,20){\line(0,1){20}}
\end{picture}
}
@>>>
\raisebox{-20pt}{
\begin{picture}(50,40)
\fboxsep=0mm
\put(0,30){\colorbox[gray]{0.7}{\makebox(10,10){}}}
\put(10,30){\colorbox[gray]{0.7}{\makebox(10,10){}}}
\put(30,20){\colorbox[gray]{0.7}{\makebox(10,10){}}}
\put(10,20){\colorbox[gray]{0.7}{\makebox(10,10){}}}
\put(0,40){\line(1,0){50}}
\put(0,30){\line(1,0){50}}
\put(10,20){\line(1,0){40}}
\put(20,10){\line(1,0){20}}
\put(30,0){\line(1,0){10}}
\put(0,30){\line(0,1){10}}
\put(10,20){\line(0,1){20}}
\put(20,10){\line(0,1){30}}
\put(30,0){\line(0,1){40}}
\put(40,0){\line(0,1){40}}
\put(50,20){\line(0,1){20}}
\put(20,30){\makebox(10,10){$\times$}}
\end{picture}
}
\\[5pt]
@VVV
@VVV
\\[5pt]
\raisebox{-20pt}{
\begin{picture}(50,40)
\fboxsep=0mm
\put(0,30){\colorbox[gray]{0.7}{\makebox(10,10){}}}
\put(10,30){\colorbox[gray]{0.7}{\makebox(10,10){}}}
\put(20,30){\colorbox[gray]{0.7}{\makebox(10,10){}}}
\put(20,10){\colorbox[gray]{0.7}{\makebox(10,10){}}}
\put(0,40){\line(1,0){50}}
\put(0,30){\line(1,0){50}}
\put(10,20){\line(1,0){40}}
\put(20,10){\line(1,0){20}}
\put(30,0){\line(1,0){10}}
\put(0,30){\line(0,1){10}}
\put(10,20){\line(0,1){20}}
\put(20,10){\line(0,1){30}}
\put(30,0){\line(0,1){40}}
\put(40,0){\line(0,1){40}}
\put(50,20){\line(0,1){20}}
\put(10,20){\makebox(10,10){$\times$}}
\end{picture}
}
@>>>
\raisebox{-20pt}{
\begin{picture}(50,40)
\fboxsep=0mm
\put(0,30){\colorbox[gray]{0.7}{\makebox(10,10){}}}
\put(10,30){\colorbox[gray]{0.7}{\makebox(10,10){}}}
\put(30,20){\colorbox[gray]{0.7}{\makebox(10,10){}}}
\put(20,10){\colorbox[gray]{0.7}{\makebox(10,10){}}}
\put(0,40){\line(1,0){50}}
\put(0,30){\line(1,0){50}}
\put(10,20){\line(1,0){40}}
\put(20,10){\line(1,0){20}}
\put(30,0){\line(1,0){10}}
\put(0,30){\line(0,1){10}}
\put(10,20){\line(0,1){20}}
\put(20,10){\line(0,1){30}}
\put(30,0){\line(0,1){40}}
\put(40,0){\line(0,1){40}}
\put(50,20){\line(0,1){20}}
\put(20,30){\makebox(10,10){$\times$}}
\put(10,20){\makebox(10,10){$\times$}}
\end{picture}
}
@>>>
\raisebox{-20pt}{
\begin{picture}(50,40)
\fboxsep=0mm
\put(0,30){\colorbox[gray]{0.7}{\makebox(10,10){}}}
\put(20,20){\colorbox[gray]{0.7}{\makebox(10,10){}}}
\put(20,10){\colorbox[gray]{0.7}{\makebox(10,10){}}}
\put(30,20){\colorbox[gray]{0.7}{\makebox(10,10){}}}
\put(0,40){\line(1,0){50}}
\put(0,30){\line(1,0){50}}
\put(10,20){\line(1,0){40}}
\put(20,10){\line(1,0){20}}
\put(30,0){\line(1,0){10}}
\put(0,30){\line(0,1){10}}
\put(10,20){\line(0,1){20}}
\put(20,10){\line(0,1){30}}
\put(30,0){\line(0,1){40}}
\put(40,0){\line(0,1){40}}
\put(50,20){\line(0,1){20}}
\put(10,30){\makebox(10,10){$\times$}}
\end{picture}
}
@>>>
\raisebox{-20pt}{
\begin{picture}(50,40)
\fboxsep=0mm
\put(10,20){\colorbox[gray]{0.7}{\makebox(10,10){}}}
\put(20,20){\colorbox[gray]{0.7}{\makebox(10,10){}}}
\put(30,20){\colorbox[gray]{0.7}{\makebox(10,10){}}}
\put(20,10){\colorbox[gray]{0.7}{\makebox(10,10){}}}
\put(0,40){\line(1,0){50}}
\put(0,30){\line(1,0){50}}
\put(10,20){\line(1,0){40}}
\put(20,10){\line(1,0){20}}
\put(30,0){\line(1,0){10}}
\put(0,30){\line(0,1){10}}
\put(10,20){\line(0,1){20}}
\put(20,10){\line(0,1){30}}
\put(30,0){\line(0,1){40}}
\put(40,0){\line(0,1){40}}
\put(50,20){\line(0,1){20}}
\put(0,30){\makebox(10,10){$\times$}}
\end{picture}
}
\\[5pt]
@VVV
@VVV
@VVV
@VVV
\\[5pt]
\raisebox{-20pt}{
\begin{picture}(50,40)
\fboxsep=0mm
\put(0,30){\colorbox[gray]{0.7}{\makebox(10,10){}}}
\put(10,30){\colorbox[gray]{0.7}{\makebox(10,10){}}}
\put(20,30){\colorbox[gray]{0.7}{\makebox(10,10){}}}
\put(30,0){\colorbox[gray]{0.7}{\makebox(10,10){}}}
\put(0,40){\line(1,0){50}}
\put(0,30){\line(1,0){50}}
\put(10,20){\line(1,0){40}}
\put(20,10){\line(1,0){20}}
\put(30,0){\line(1,0){10}}
\put(0,30){\line(0,1){10}}
\put(10,20){\line(0,1){20}}
\put(20,10){\line(0,1){30}}
\put(30,0){\line(0,1){40}}
\put(40,0){\line(0,1){40}}
\put(50,20){\line(0,1){20}}
\put(10,20){\makebox(10,10){$\times$}}
\put(20,10){\makebox(10,10){$\times$}}
\end{picture}
}
@>>>
\raisebox{-20pt}{
\begin{picture}(50,40)
\fboxsep=0mm
\put(0,30){\colorbox[gray]{0.7}{\makebox(10,10){}}}
\put(10,30){\colorbox[gray]{0.7}{\makebox(10,10){}}}
\put(30,20){\colorbox[gray]{0.7}{\makebox(10,10){}}}
\put(30,0){\colorbox[gray]{0.7}{\makebox(10,10){}}}
\put(0,40){\line(1,0){50}}
\put(0,30){\line(1,0){50}}
\put(10,20){\line(1,0){40}}
\put(20,10){\line(1,0){20}}
\put(30,0){\line(1,0){10}}
\put(0,30){\line(0,1){10}}
\put(10,20){\line(0,1){20}}
\put(20,10){\line(0,1){30}}
\put(30,0){\line(0,1){40}}
\put(40,0){\line(0,1){40}}
\put(50,20){\line(0,1){20}}
\put(20,30){\makebox(10,10){$\times$}}
\put(10,20){\makebox(10,10){$\times$}}
\put(20,10){\makebox(10,10){$\times$}}
\end{picture}
}
@>>>
\raisebox{-20pt}{
\begin{picture}(50,40)
\fboxsep=0mm
\put(0,30){\colorbox[gray]{0.7}{\makebox(10,10){}}}
\put(20,20){\colorbox[gray]{0.7}{\makebox(10,10){}}}
\put(30,20){\colorbox[gray]{0.7}{\makebox(10,10){}}}
\put(30,0){\colorbox[gray]{0.7}{\makebox(10,10){}}}
\put(0,40){\line(1,0){50}}
\put(0,30){\line(1,0){50}}
\put(10,20){\line(1,0){40}}
\put(20,10){\line(1,0){20}}
\put(30,0){\line(1,0){10}}
\put(0,30){\line(0,1){10}}
\put(10,20){\line(0,1){20}}
\put(20,10){\line(0,1){30}}
\put(30,0){\line(0,1){40}}
\put(40,0){\line(0,1){40}}
\put(50,20){\line(0,1){20}}
\put(10,30){\makebox(10,10){$\times$}}
\put(20,10){\makebox(10,10){$\times$}}
\end{picture}
}
@>>>
\raisebox{-20pt}{
\begin{picture}(50,40)
\fboxsep=0mm
\put(10,20){\colorbox[gray]{0.7}{\makebox(10,10){}}}
\put(20,20){\colorbox[gray]{0.7}{\makebox(10,10){}}}
\put(30,20){\colorbox[gray]{0.7}{\makebox(10,10){}}}
\put(30,0){\colorbox[gray]{0.7}{\makebox(10,10){}}}
\put(0,40){\line(1,0){50}}
\put(0,30){\line(1,0){50}}
\put(10,20){\line(1,0){40}}
\put(20,10){\line(1,0){20}}
\put(30,0){\line(1,0){10}}
\put(0,30){\line(0,1){10}}
\put(10,20){\line(0,1){20}}
\put(20,10){\line(0,1){30}}
\put(30,0){\line(0,1){40}}
\put(40,0){\line(0,1){40}}
\put(50,20){\line(0,1){20}}
\put(0,30){\makebox(10,10){$\times$}}
\put(20,10){\makebox(10,10){$\times$}}
\end{picture}
}
\end{CD}
$$
\caption{Excited diagrams in $S(5,4,2,1)$ viewed as a type $B$ heap}
\label{fig:shifted-B_excited}
\end{figure}
\end{exam}
\section{%
Equivariant $K$-theory of Kac--Moody partial flag varieties
}

In this section we review the basic properties of the equivariant $K$-theory 
of thick flag varieties following \cite[Section~3]{LSS}, 
and rephrase the Billey-type formula and the Chevalley-type formula 
in terms of combinatorics of $d$-complete posets.

\subsection{
Equivariant $K$-theory and localization
}

Let $A = (a_{ij})_{i,j \in I}$ be a symmetrizable generalized Cartan matrix, 
and $\Gamma$ the corresponding Dynkin diagram with node set $I$.
Then the associated Kac--Moody group $\GG$ over $\Comp$ is constructed from the following data:
the weight lattice ($\Int$-module) $\Lambda$, the simple roots $\Pi = \{ \alpha_i : i \in I \}$, 
the simple coroots $\Pi^\vee = \{ \alpha^\vee_i :i \in I \}$, 
and the fundamental weights $\{ \lambda_i : i \in I \}$ (see the beginning of Subsection~2.2).

In what follows, we fix a subset $J$ of $I$.
Let $\BB$ be a Borel subgroup corresponding to the positive system $\Phi_+$ 
and $\TT \subset \BB$ a maximal torus.
Let $\PP_-$ be the opposite parabolic subgroup corresponding to the subset $J$, 
which contains the opposite Borel subgroup $\BB_-$ such that $\BB \cap \BB_- = \TT$.
Then we can introduce the Kashiwara thick partial flag variety $\XX = \GG/\PP_-$.
(We refer the readers to \cite{Ka} for a construction of $\XX$.)

Let $W_J$ be the parabolic subgroup of $W$ corresponding to $J$ 
and $W^J$ be the set of minimum length coset representatives of $W/W_J$.
For each element $v \in W^J$, we put $\XX^\circ_v = \BB v \PP_-/\PP_-$ 
and $\XX_v = \overline{\XX^\circ_v}$, the Zariski closure of $\XX^\circ_v$, 
which are called the Schubert cell and the Schubert variety respectively.
Then $\XX_v$ has codimension $l(v)$ in $\XX$ and
$$
\XX_v = \bigsqcup_{w \in W^J, \, w \ge v} \XX^\circ_w.
$$

Let $K_\TT(\XX)$ be the $\TT$-equivariant $K$-theory of $\XX$.
Then $K_\TT(\XX)$ has a commutative associative $K_\TT(\point)$-algebra structure.
Here the $\TT$-equivalent $K$-theory $K_\TT(\point)$ of a point is isomorphic 
to the group algebra $\Int[\Lambda]$ with basis $\{ e^\lambda : \lambda \in \Lambda \}$,  
and to the representation ring $R(\TT)$ of $\TT$.
For each $v \in W^J$, let $[\OO_v]$ be the class of the structure sheaf $\OO_v$ of $\XX_v$ 
in $K_\TT(\XX)$ and call it the equivariant Schubert class.
Then we have
$$
K_\TT(\XX) \cong \prod_{v \in W^J} K_\TT(\point) [\OO_v].
$$
Any elements of $K_\TT(\XX)$ is a (possibly infinite) $K_\TT(\point)$-linear combination 
of the equivariant Schubert classes.

Each $w \in W^J$ gives a $\TT$-fixed point $e_w = w \PP_-/\PP_- \in \XX$.
Then the inclusion map $\iota_w : \{ e_w \} \to \XX$ induces the pull-back ring homomorphism, 
called the localization map at $w$,
$$
\iota_w^* : K_{\TT}(\XX) \to K_{\TT}(e_w) \cong \Int[\Lambda].
$$
If $\LL^\lambda$ is the line bundle on $\XX$ corresponding to a weight $\lambda \in \Lambda$, 
then the image of the class $[\LL^\lambda]$ under the localization map is given by
$\iota_w^*( [\LL^\lambda] ) = e^{w \lambda}$.
For two elements $v$, $w \in W^J$, we denote by $\xi^v|_w$ 
the image of the $\TT$-equivariant Schubert class $\xi^v = [\OO_v] \in K_{\TT}(\XX)$ 
under the localization map $\iota_w^*$:
$$
\xi^v|_w = \iota_w^* ( [\OO_v] ).
$$
Then the Billey-type formula for the equivariant $K$-theory can be stated as follows:

\begin{prop}
\label{prop:Billey}
\textup{(\cite[Proposition~2.10]{LSS})}
Let $v$, $w \in W^J$, and fix a reduced expression $w = s_{i_1} s_{i_2} \dots s_{i_N}$ of $w$.
Then we have
\begin{equation}
\label{eq:Billey}
\xi^v|_w 
 =
\sum_{(k_1, \dots, k_r)}
 (-1)^{r-l(v)} \prod_{a=1}^r \left( 1 - e^{\beta^{(k_a)}} \right),
\end{equation}
where the summation is taken over all sequences $(k_1, \dots, k_r)$ 
such that $1 \le k_1 < k_2 < \dots < k_r \le N$ and $s_{i_{k_1}} * \dots * s_{i_{k_r}} = v$ 
\textup{(}with respect to the Demazure product\textup{)}, 
and $\beta^{(k)}$ is given by $\beta^{(k)} = s_{i_1} \dots s_{i_{k-1}}(\alpha_{i_k})$ for $1 \le k \le N$.
\end{prop}

By using Lemma~\ref{lem:red} (a), we can deduce the following corollary from Proposition~\ref{prop:Billey}.

\begin{coro}
\label{cor:Billey}
\begin{enumerate}
\item[(a)]
For $w \in W^J$, we have
\begin{equation}
\label{eq:Billey2}
\xi^w|_w = \prod_{k=1}^N \left( 1 - e^{\beta^{(k)}} \right).
\end{equation}
In particular, $\xi^w|_w \neq 0$.
\item[(b)]
Let $v$, $w \in W^J$.
If $\xi^v|_w \neq 0$, then we have $v \le w$ in the Bruhat order.
\end{enumerate}
\end{coro}

\subsection{%
Equivariant $K$-theoretical Littlewood--Richardson coefficients
}

We consider the structure constants for the multiplication in $K_{\TT}(\XX)$ 
with respect to the equivariant Schubert classes.
For $u$, $v$, $w \in W^J$, we denote by $c^w_{u,v} \in K_{\TT}(\point)$ 
the structure constant determined by
$$
[\OO_u] [\OO_v] = \sum_{w \in W^J} c^w_{u,v} [\OO_w].
$$

\begin{lemm}
\label{lem:LR}
If $c^w_{u,v} \neq 0$, then $u \le w$ and $v \le w$.
\end{lemm}

\begin{proof}
We use the induction on $l(w)$ to prove that, 
if $u \not\le w$ or $v \not\le w$, then $c^w_{u,v} = 0$.
Assume that $u \not\le w$ or $v \not\le w$.
By apply the localization map $\iota_w^*$ to 
$[\OO_u] [\OO_v] = \sum_{x \in W^J} c^x_{u,v} [\OO_x]$ 
and then by using Corollary~\ref{cor:Billey} (b), we have
$$
\left( \xi^u|_w \right) \cdot \left( \xi^v|_w \right)
 =
\sum_{x \le w} c^x_{u,v} \xi^x|_w.
$$
If there exists an element $x \in W^J$ satisfying $x < w$ and $c^x_{u,v} \neq 0$, 
then we have $u \le x$ and $v \le x$ by the induction hypothesis, 
and hence $u \le w$ and $v \le w$, which contradicts to the assumption.
Hence, by using Corollary~\ref{cor:Billey} (b) and the assumption, we have
$$
0 = c^w_{u,v} \xi^w|_w.
$$
Since $\xi^w_w \neq 0$ (Corollary~\ref{cor:Billey} (a)), we obtain $c^w_{u,v} = 0$.
\end{proof}

\begin{prop}
\label{prop:LR=xi}
For $v$, $w \in W^J$, we have
\begin{equation}
\label{eq:LR=xi}
c^w_{v,w} = \xi^v|_w.
\end{equation}
\end{prop}

\begin{proof}
By apply the localization map $\iota_w^*$ to 
$[\OO_v] [\OO_w] = \sum_{x \in W^J} c^x_{v,w} [\OO_x]$ 
and then by using Corollary~\ref{cor:Billey} (b), we have
$$
\left( \xi^v|_w \right) \cdot \left( \xi^w|_w \right)
 =
\sum_{x \le w} c^x_{v,w} \xi^x|_w.
$$
By Lemma~\ref{lem:LR}, we see that $c^x_{v,w} = 0$ unless $w \le x$.
By Corollary~\ref{cor:Billey} (b), we see that $\xi^x|_w = 0$ unless $x \le w$.
Hence we have
$$
\left( \xi^v|_w \right) \cdot \left( \xi^w|_w \right)
 =
c^w_{v,w} \xi^w|_w.
$$
Since $\xi^w|_w \neq 0$ (Corollary~\ref{cor:Billey} (a)), we obtain the desired equality.
\end{proof}

The following lemma gives a recurrence of the equivariant $K$-theoretical 
Littlewood--Richardson coefficients $c^w_{u,v}$.
We use the same idea as \cite[Corollary~6.5]{M} and \cite[Proposition~3.1]{PY}.

\begin{lemm}
\label{lem:recLR}
Let $u$, $v$, $w \in W^J$ and $s \in W^J$ a simple reflection.
If $c^w_{s,w} \neq c^u_{s,u}$, then we have
$$
c^w_{u,v}
 =
\frac{ 1 }
     { c^w_{s,w} - c^u_{s,u} }
\left(
 \sum_{u < x \le w} c^x_{s,u} c^w_{x,v}
 -
 \sum_{u \le y < w} c^w_{s,y} c^y_{u,v}
\right).
$$
In particular, we have
\begin{equation}
\label{eq:recLR}
c^w_{u,w}
 =
\frac{ 1 }
     { c^w_{s,w} - c^u_{s,u} }
\sum_{u < x \le w}
 c^x_{s,u} c^w_{x,w}.
\end{equation}
\end{lemm}

\begin{proof}
Consider the associativity
$$
\left( [\OO_s] [\OO_u] \right) [\OO_v]
 =
[\OO_s] \left( [\OO_u] [\OO_v] \right).
$$
Taking the coefficients of $[\OO_w]$ in the both hand sides and using Lemma~\ref{lem:LR}, 
we have
$$
c^u_{s,u} c^w_{u,v}
 +
\sum_{u < x \le w} c^x_{s,u} c^w_{x,v}
 =
c^w_{s,w} c^w_{u,v}
 +
\sum_{u \le y < w} c^w_{s,y} c^y_{u,v},
$$
from which we get the conclusion. 
\end{proof}

The Chevalley formula gives a combinatorial expression of $c^w_{s,v}$ for a simple reflection $s$.
To state the Chevalley formula of \cite{LS} we need several notations.
For a dominant weight $\lambda \in \Lambda$, we put
$$
\mathbb{H}_\lambda
 =
\{ (\gamma^\vee, k) : 
 \gamma^\vee \in \Phi^\vee_+, \ k \in \Int, \ 0  \le k < \langle \gamma^\vee, \lambda \rangle
\}.
$$
Fix a total order on $I$ so that $I = \{ i_1 < \cdots < i_r \}$, 
and define a map $\iota : \mathbb{H}_\lambda \to \Rat^{r+1}$ by
$$
\iota \left(\gamma^\vee , k \right)
 =
\frac{ 1 }
     { \langle \gamma^\vee, \lambda \rangle }
\left( k, c_1, \cdots, c_r \right).
$$
where $c_1, \dots, c_r$ are the coefficients of the simple roots in $\gamma^\vee$ 
given by $\gamma^\vee = \sum_{j=1}^r c_j \alpha^\vee_{i_j}$.
Then it is known that $\iota$ is injective.
We define a total ordering $<$ on $\mathbb{H}_\lambda$ by
$$
h < h' \Longleftrightarrow \iota(h) <_{\text{lex}} \iota(h'),
$$
where $<_{\text{lex}}$ is the lexicographical ordering on $\Rat^{r+1}$.
For $h = (\gamma^\vee, k)$, we define affine transformations $r_h$ and $\tilde{r}_h$ on $\Lambda$ by
\begin{align*}
r_h(\mu) &= \mu - \langle \gamma^\vee, \mu \rangle \gamma,
\\
\tilde{r}_h(\mu) &= r_h(\mu) + \left( \langle \gamma^\vee, \lambda \rangle - k \right) \gamma.
\end{align*}
Note that $r_h = s_{\gamma}$.
Now we can state the Chevalley formula for the equivariant $K$-theory of the partial flag variety $\XX$.

\begin{prop}
\label{prop:Chevalley}
\textup{(\cite[Theorem~4.8 (4.12) and (4.13)]{LS}, see also \cite[Corollary~7.1]{LP})}
Let $s$ be a simple reflection such that $s \in W^J$ and $v$, $w \in W^J$.
If $s = s_i$ and $\lambda = \lambda_i$ is the corresponding fundamental weight, then we have
\begin{equation}
\label{eq:Chevalley}
c^w_{s,v}
 =
\begin{cases}
 1 - e^{\lambda - v \lambda} &\text{if $w=v$},
\\
 \displaystyle\sum_{(h_1, \cdots, h_r)}
  (-1)^{r-1}
  e^{\lambda - v \tilde{r}_{h_1} \cdots \tilde{r}_{h_r} \lambda}
 &\text{if $w > v$,}
\\
 0 &\text{otherwise,}
\end{cases}
\end{equation}
where the summation is taken over all sequences $(h_1, \cdots, h_r)$ of length $r \ge 1$ 
satisfying the following two conditions:
\begin{enumerate}
\item[(H1)]
$h_1 > h_2 > \cdots > h_r$ in $\mathbb{H}_\lambda$,
\item[(H2)]
$v \lessdot v r_{h_1} \lessdot v r_{h_1} r_{h_2} \lessdot \cdots \lessdot v r_{h_1} \cdots r_{h_r} = w$ 
is a saturated chain in $W^J$.
\end{enumerate}
\end{prop}

\subsection{%
Connection to $d$-complete posets
}

In this subsection we rephrase the Billey-type formula and the Chevalley-type formula 
in terms of combinatorics of $d$-complete posets.

Let $P$ be a connected $d$-complete poset with top tree $\Gamma$.
We regard $\Gamma$ as a simply-laced Dynkin diagram with node set $I$.
Let $\alpha_P$ and $\lambda_P$ be the simple root and the fundamental weight 
corresponding to the color $i_P$ of the maximum element of $P$.
We apply the results of Subsections~4.1 and 4.2 
to the Kashiwara thick partial flag variety 
$\XX = \GG/\PP_-$, where $\PP_-$ is the maximal parabolic subgroup corresponding 
to $J = I \setminus \{ i_P \}$.
In this case, the parabolic subgroup $W_J$ coincides with the stabilizer $W_{\lambda_P}$ of $\lambda_P$ in $W$, 
and the minimum length coset representatives $W^J$ is denoted by $W^{\lambda_P}$.

By using a labeling of the elements of $P$ with $p_1, \cdots, p_N$ ($N= \# P$) 
so that $p_i < p_j$ in $P$ implies $i<j$, 
we can associate to each subset $D = \{ i_1, \dots, i_r \}$ ($i_1 < \dots < i_r$) of $P$ 
a well-defined element $w_D = s(p_{i_1}) \cdots s(p_{i_r}) \in W$.
Then the following formula is obtained from the Billey-type formula.

\begin{prop}
\label{prop:Billey_P}
Let $P$ be a connected $d$-complete poset and $F$ an order filter of $P$.
Then we have
\begin{equation}
\label{eq:Billey_P}
\xi^{w_F}|_{w_P}
 =
\sum_{E \in \EE^*_P(F)}
 (-1)^{\# E - \# F}
 \prod_{p \in E} \left( 1 - \vectz[H_P(p)] \right),
\end{equation}
under the identification $z_i = e^{\alpha_i}$ \textup{(}$i \in I$\textup{)}.
\end{prop}

\begin{proof}
Follows from Proposition~\ref{prop:Billey} 
by using Proposition~\ref{prop:root} (b) and Proposition~\ref{prop:E*=R*}.
\end{proof}

Also the following explicit expression is obtained from the Chevalley-type formula.

\begin{prop}
\label{prop:Chevalley_P}
Let $P$ be a connected $d$-complete poset and put $s = s_{i_P}$.
For two order filters $F$ and $F'$ of $P$, we have
\begin{equation}
\label{eq:Chevalley_P}
c^{w_{F'}}_{s,w_F}
 =
\begin{cases}
 1 - \vectz[F]
 &\text{if $F' = F$,}
\\
 (-1)^{\# (F' \setminus F) -1} \vectz[F]
 &\text{if $F' \supsetneq F$ and $F' \setminus F$ is an antichain,}
\\
 0
 &\text{otherwise,}
\end{cases}
\end{equation}
under the identification $z_i = e^{\alpha_i}$ \textup{(}$i \in I$\textup{)}.
\end{prop}

First we consider the case $r=1$ in Proposition~\ref{prop:Chevalley}.

\begin{lemm}
\label{lem:r=1}
Let $F$ be an order filter of $P$ and $h = (\gamma^\vee,k) \in \mathbb{H}_{\lambda_P}$.
If $w_F r_h \in W^{\lambda_P}$ and $w_F \lessdot w_F r_h \le w_P$, 
then there exists $p \in P$ such that $F' = F \sqcup \{ p \}$ is an order filter of $P$,  
$w_F r_h = w_{F'}$ and $\gamma^\vee = \gamma^\vee(p)$.
In this case $k=0$ and $\tilde{r}_h \lambda_P = \lambda_P$.
\end{lemm}

\begin{proof}
Since the interval $[e, w_P]$ in $W^{\lambda_P}$ is isomorphic to the poset of order filters of $P$ 
(Proposition~\ref{prop:Weyl} (a)), 
there exists a unique order filter $F'$ of $P$ such that $F' \supset F$, 
$\# F' = \# F + 1$ and $w_{F'} = w_F r_h$.
Hence we have $p \in P$ such that $F' = F \sqcup \{ p \}$ and $w_{F'} = s(p) w_F$.
We take a linear extension of $P$ such that $F = \{ p_{n+1}, \cdots, p_N \}$ 
with $N = \# P$ and $n = \# (P \setminus F)$.
If $p = p_m$, then $p$ is incomparable with $p_{m+1}, \cdots, p_n$, 
hence $s(p_m)$ is commutative with $s(p_{m+1}), \cdots, s(p_n)$ by Property (C5) 
in Proposition~\ref{prop:dc_color}, 
and thus $s(p_i)\alpha^\vee(p_m)=\alpha^\vee(p_m)$ for $m+1 \le i \le n$.
Hence we have
\begin{align*}
\gamma^\vee
 &= 
w_F^{-1} \alpha^\vee(p_m)
\\
 &=
s(p_N) \cdots s(p_{n+1}) \alpha^\vee(p_m)
\\
 &=
s(p_N) \cdots s(p_{n+1}) s(p_n) \cdots s(p_{m+1}) \alpha^\vee(p_m)
\\
 &=
\gamma^\vee(p_m).
\end{align*}
By Proposition~\ref{prop:root} (c), we see that $k=0$ and $\tilde{r}_h \lambda_P = \lambda_P$.
\end{proof}

Now we deduce Proposition~\ref{prop:Chevalley_P} from the Chevalley-type formula.

\begin{proof}[Proof of Proposition~\ref{prop:Chevalley_P}]
It follows from Proposition~\ref{prop:Weyl} (b)
and Proposition~\ref{prop:Chevalley} that
$$
c^{w_F}_{s, w_F} = 1 - \vectz[F].
$$

Suppose that there exists a sequence $(h_1, \cdots, h_r)$ of elements in $\mathbb{H}_{\lambda_P}$ 
satisfying Conditions (H1) and (H2) in Proposition~\ref{prop:Chevalley}.
Then by Lemma~\ref{lem:r=1}, we have a sequence $(q_1, \cdots, q_r)$ of elements of $P$ such that 
$F_i = F \sqcup \{ q_1, \cdots, q_i \}$ is an order filter of $P$, $h_i = (\gamma^\vee(q_i),0)$ 
for $1 \le i \le r$ and $\tilde{r}_{h_1} \cdots \tilde{r}_{h_r} \lambda_P = \lambda_P$.
Now we show that $\{ q_1, \cdots, q_r \}$ is an antichain.
Assume to the contrary that there exist $i$ and $j$ such that $i<j$ and $q_i$ and $q_j$ are comparable.
Since $q_i$ is maximal in $P \setminus (F \sqcup \{ q_1, \cdots, q_{i-1} \})$ 
and $q_j \in P \setminus (F \sqcup \{ q_1, \cdots, q_{i-1} \})$, 
we have $q_i > q_j$.
Then by Proposition~\ref{prop:root} (a), we see that $\gamma^\vee(q_i) < \gamma^\vee(q_j)$.
Hence by the definition of the total order on $\mathbb{H}_{\lambda_P}$, we have $h_i < h_j$,
which contradicts to Condition (H1).
Moreover it follows from Proposition~\ref{prop:Weyl} (b) that
$$
e^{\lambda_P - w_F \tilde{r}_{h_1} \cdots \tilde{r}_{h_r} \lambda_P}
 = \vectz[F'].
$$

Conversely, suppose that $F' \supsetneq F$ and $F' \setminus F$ is an antichain.
For $q \in F' \setminus F$, we put $h(q) = (\gamma^\vee(q),0) \in \mathbb{H}_{\lambda_P}$.
Since $F' \setminus F$ is an antichain, we can label the elements of $F' \setminus F$ 
so that $h(q_1) > \cdots > h(q_r)$.
Then $(h(q_1), \cdots, h(q_r))$ is the unique sequence satisfying Conditions (H1) and (H2) 
in Proposition~\ref{prop:Chevalley}.
\end{proof}

\section{%
Proof and corollaries of Main Theorem
}

In this section, we give a proof of the main theorem (Theorem~\ref{thm:main} in Introduction) 
and derive several consequences.

\subsection{
Proof of the Main Theorem
}

Recall the main theorem of this paper:

\begin{theo}
\label{thm:main1}
Let $P$ be a $d$-complete poset and $F$ an order filter.
Then the multivariate generating function of $(P \setminus F)$-partitions, 
where $P \setminus F$ is viewed as an induced subposet of $P$, is given by
\begin{equation}
\label{eq:main1}
\sum_{\sigma \in \AP(P \setminus F)} \vectz^\sigma
 =
\sum_{D \in \EE_P(F)}
 \frac{ \prod_{v \in B(D)} \vectz[H_P(v)] }
      { \prod_{v \in P \setminus D} ( 1 - \vectz[H_P(v)] ) },
\end{equation}
where $D$ runs over all excited diagrams of $F$ in $P$.
\end{theo}

Since the both sides of \eqref{eq:main1} factor into the product over the connected components of $P$, 
the assertion of this theorem follows easily from the case where $P$ is a connected $d$-complete poset.
Hence Theorem~\ref{thm:main1} is a direct consequence of the following two theorems, 
which describe the ratio $\xi^{w_F}|_{w_P} \big/ \xi^{w_P}|_{w_P}$ of the localizations of elements 
in the equivariant $K$-theory $K_\TT(\XX)$ in two ways.

\begin{theo}
\label{thm:xi1}
For a connected $d$-complete poset $P$ and an order filter $F$ of $P$, we have
\begin{equation}
\label{eq:xi1}
\frac{ \xi^{w_F}|_{w_P} }
     { \xi^{w_P}|_{w_P} }
 =
\sum_{\sigma \in \AP(P \setminus F)} \vectz^\sigma,
\end{equation}
under the identification $z_i = e^{\alpha_i}$ \textup{(}$i \in I$\textup{)}.
\end{theo}

\begin{theo}
\label{thm:xi2}
For a connected $d$-complete poset $P$ and an order filter $F$ of $P$, we have
\begin{equation}
\label{eq:xi2}
\frac{ \xi^{w_F}|_{w_P} }
     { \xi^{w_P}|_{w_P} }
 =
\sum_{D \in \EE_P(F)}
 \frac{ \prod_{q \in B(D)} \vectz[H_P(q)] }
      { \prod_{p \in P \setminus D} (1 - \vectz[H_P(p)]) },
\end{equation}
under the identification $z_i = e^{\alpha_i}$ \textup{(}$i \in I$\textup{)}.
\end{theo}

First we prove Theorem~\ref{thm:xi1} by using the Chevalley-type formula (Proposition~\ref{prop:Chevalley_P}).

\begin{proof}[Proof of Theorem~\ref{thm:xi1}]
For an order filter $F$ of $P$, we put
$$
Z_{P/F}(\vectz)
 = 
\frac{ \xi^{w_F}|_{w_P} }
     { \xi^{w_P}|_{w_P} },
\quad
G_{P/F}(\vectz) = \sum_{\sigma \in \AP(P \setminus F)} \vectz^\sigma.
$$
It is clear that $Z_{P/P}(\vectz) = G_{P/P}(\vectz) = 1$.
Hence it is enough to show that $Z_{P/F}(\vectz)$ and $G_{P/F}(\vectz)$ satisfy the same recurrences:
\begin{align}
\label{eq:recZ}
Z_{P/F}(\vectz)
 &=
\frac{ 1 }
     { 1 - \vectz[P \setminus F] }
\sum_{F'}
 (-1)^{\# (F' \setminus F) -1 }
 Z_{P/F'}(\vectz),
\\
\label{eq:recG}
G_{P/F}(\vectz)
 &=
\frac{ 1 }
     { 1 - \vectz[P \setminus F] }
\sum_{F'}
 (-1)^{\# (F' \setminus F) -1 }
 G_{P/F'}(\vectz),
\end{align}
where $F'$ runs over all order filters such that $F \subsetneq F' \subset P$ and 
$F' \setminus F$ is an antichain.

First we prove \eqref{eq:recZ}.
Under the isomorphism of posets given in Proposition~\ref{prop:Weyl} (a), 
the interval $(w_F,w_P] = \{ z \in W^{\lambda_P} : w_F < z \le w_P \}$ corresponds 
to $\{ F' : \text{$F'$ is an order filter of $P$ and $F \subsetneq F' \subset P$} \}$.
Then by using the recurrence \eqref{eq:recLR} and Proposition~\ref{prop:Chevalley_P}, we see that
\begin{align*}
\xi^{w_F}|_{w_P}
 &=
\frac{ 1 }
     { (1 - \vectz[P]) - (1 - \vectz[F]) }
\sum_{F'}
 (-1)^{\# (F' \setminus F) - 1} \vectz[F] \xi^{w_{F'}}|_{w_P}
\\
 &=
\frac{ 1 }
     { 1 - \vectz[P \setminus F] }
\sum_{F'}
 (-1)^{\# (F' \setminus F) - 1} \xi^{w_{F'}}|_{w_P}.
\end{align*}

Next we prove \eqref{eq:recG}.
Let $M$ be the set of maximal elements of $P \setminus F$.
Then we have
$$
\sum_{F'}
 (-1)^{\# (F' \setminus F) -1 }
 G_{P/F'}(\vectz)
 =
\sum_{I \subset M, I \neq \emptyset}
 (-1)^{\# I - 1}
 G_{P/(F \sqcup I)}(\vectz).
$$
For $I \subset M$, we put
$$
\AP(P \setminus F)_I 
 =
\{ \sigma \in \AP(P \setminus F) : \sigma(x) = 0 \text{ for all $x \in I$} \}.
$$
Then we have
$$
G_{P/(F \sqcup I)}(\vectz) = \sum_{\sigma \in \AP(P \setminus F)_I} \vectz^\sigma.
$$
By the Inclusion-Exclusion Principle, we have
$$
\sum_{F'}
 (-1)^{\# (F' \setminus F) -1 }
 G_{P/F'}(\vectz)
 =
\sum_{\sigma \in \AP'(P \setminus F)} \vectz^\sigma,
$$
where we put
$$
\AP'(P \setminus F)
 =
\{ \sigma \in \AP(P \setminus F) :
 \text{$\sigma(x) = 0$ for some $x \in M$}
\}.
$$
Given $\sigma \in \AP(P \setminus F)$, let $m = \min \{ \sigma(x) : x \in P \setminus F \}$ 
and define $\sigma' \in \AP(P \setminus F)$ by $\sigma'(x) = \sigma(x) - m$ ($x \in P \setminus F$).
Then the map $\sigma \mapsto (m, \sigma')$ gives a bijection from $\AP(P \setminus F)$ 
to $\Nat \times \AP'(P \setminus F)$, and
$$
\vectz^\sigma = \vectz[P \setminus F]^m \cdot \vectz^{\sigma'}.
$$
Hence we have
$$
\sum_{\sigma \in \AP(P \setminus F)} \vectz^\sigma
  =
\frac{ 1 }
     { 1 - \vectz[P \setminus F] }
\sum_{\sigma \in \AP'(P \setminus F)} \vectz^\sigma.
$$
This completes the proof.
\end{proof}

Next we derive Theorem~\ref{thm:xi2} from the Billey-type formula (Proposition~\ref{prop:Billey_P}).

\begin{proof}[Proof of Theorem~\ref{thm:xi2}]
By Proposition~\ref{prop:Billey_P}, we have
$$
\xi^{w_F}|_{w_P}
 =
\sum_{E \in \EE^*_P(F)}
 (-1)^{\# E - \# F}
 \prod_{p \in E} \left( 1 - \vectz[H_P(p)] \right).
$$
By using Proposition~\ref{prop:E*}, we see that
\begin{align*}
\xi^{w_F}|_{w_P}
 &=
\sum_{D \in \EE_P(F)}
 \prod_{p \in D} \left( 1 - \vectz[H_P(p)] \right)
\sum_{S \subset B(D)}
 (-1)^{\# S}
 \prod_{p \in S} \left( 1 - \vectz[H_P(p)] \right)
\\
 &=
\sum_{D \in \EE_P(F)}
 \prod_{p \in D} \left( 1 - \vectz[H_P(p)] \right)
 \prod_{p \in B(D)} \vectz[H_P(p)].
\end{align*}
By dividing the both sides by $\xi^{w_P}|_{w_P} = \prod_{p \in P} ( 1 - \vectz[H_P(p)] )$, 
we obtain the desired identity \eqref{eq:xi2}.
\end{proof}

\subsection{%
Corollaries of the Main Theorem
}

First we derive the equivariant cohomology version of Theorem~\ref{thm:main1}.
In addition to translating Nakada's colored hook formula \cite[Corollary~7.2]{N1} 
from the context of roots to the context of $d$-complete posets, 
the following corollary gives a skew generalization of it.

\begin{coro}
\label{cor:main2}
Let $P$ be a $d$-complete poset with $d$-complete coloring $c : P \to I$ 
and $F$ an order filter of $P$.
Let $\vecta = (a_i)_{i \in I}$ be indeterminates.
We put $a(p) = a_{c(p)}$ \textup{(}$p \in P$\textup{)} 
and define a linear polynomial $\vecta \langle H_P(u) \rangle$ 
as follows:
\begin{enumerate}
\item[(i)]
If $u$ is not the top of any $d_k$-interval, then we define
$$
\vecta \langle H_P(u) \rangle = \sum_{w \le u} a_{c(w)}.
$$
\item[(ii)]
If $u$ is the top of a $d_k$-interval $[v,u]$, then we define
$$
\vecta \langle H_P(u) \rangle
 =
\vecta \langle H_P(x) \rangle + \vecta \langle H_P(y) \rangle
- \vecta \langle H_P(v) \rangle,
$$
where $x$ and $y$ are the sides of $[v,u]$.
\end{enumerate}
Then we have
\begin{multline}
\label{eq:main2}
\sum_{(q_1, \dots, q_n)}
 \frac{ 1 }
      { a(q_1) (a(q_1)+a(q_2)) \cdots (a(q_1)+\dots+a(q_n)) }
\\
 =
\sum_{D \in \EE_P(F)}
 \prod_{v \in P \setminus D}
 \frac{ 1 }
      { \vecta \langle H_P(v) \rangle },
\end{multline}
where $n = \# (P\setminus F)$ and 
the summation is taken over all linear extensions of $P \setminus F$, 
\ie all labelings of the elements of $P \setminus F$ with $q_1, \dots, q_n$ so that 
 $q_i < q_j$ in $P \setminus F$ implies $i < j$.
\end{coro}

\begin{proof}
Any $(P \setminus F)$-partition $\sigma \in \AP(P \setminus F)$ is determined by 
a nonnegative integer $r\leq n$, an increasing sequence $i_1 < \cdots < i_r$ of positive integers 
and an increasing sequence $F \subset F_0 \subsetneq F_1 \subsetneq \cdots \subsetneq F_r = P$ 
of order filters of $P$, by the condition
$$
\sigma(x)
 =
\begin{cases}
 0 &\text{if $x \in F_0 \setminus F$,} \\
 i_k &\text{if $x \in F_k \setminus F_{k-1}$ and $1 \le k \le r$.}
\end{cases}
$$
Hence we have
\begin{align*}
&
\sum_{\sigma \in \AP(P \setminus F)} \vectz^\sigma
\\
&\quad
=
\sum_{F \subset F_0 \subsetneq F_1 \subsetneq \cdots \subsetneq F_r = P}
\sum_{0 < i_1 < \cdots < i_r}
 \vectz[F_1 \setminus F_0]^{i_1} \vectz[F_2 \setminus F_1]^{i_2} \cdots \vectz[P \setminus F_{r-1}]^{i_r}
\\
&\quad
=
\sum_{F \subset F_0 \subsetneq F_1 \subsetneq \cdots \subsetneq F_r = P}
 \frac{ \vectz[P \setminus F_0] }
      { 1 - \vectz[P \setminus F_0] }
 \frac{ \vectz[P \setminus F_1] }
      { 1 - \vectz[P \setminus F_1] }
 \cdots
 \frac{ \vectz[P \setminus F_{r-1}] }
      { 1 - \vectz[P \setminus F_{r-1}] }.
\end{align*}
Now by using Theorem~\ref{thm:main1}, we have
\begin{multline*}
\sum_{F \subset F_0 \subsetneq F_1 \subsetneq \cdots \subsetneq F_r = P}
 \frac{ \vectz[P \setminus F_0] }
      { 1 - \vectz[P \setminus F_0] }
 \frac{ \vectz[P \setminus F_1] }
      { 1 - \vectz[P \setminus F_1] }
 \cdots
 \frac{ \vectz[P \setminus F_{r-1}] }
      { 1 - \vectz[P \setminus F_{r-1}] }
\\
 =
\sum_{D \in \EE_P(F)}
 \frac{ \prod_{v \in B(D)} \vectz[H_P(v)] }
      { \prod_{v \in P \setminus D} ( 1 - \vectz[H_P(v)] ) }.
\end{multline*}
By substituting $z_i = t^{a_i}$ ($i \in I$) and multiplying the both sides by $(1-t)^n$, 
and then by taking the limit $t \to 1$, we obtain
\begin{multline*}
\sum_{F = F_0 \subsetneq F_1 \subsetneq \cdots \subsetneq F_n = P}
 \frac{ 1 }
      { \vecta \langle P \setminus F_0 \rangle \vecta \langle P \setminus F_1 \rangle 
        \cdots \vecta \langle P \setminus F_{n-1} \rangle
      }
\\
 =
\sum_{D \in \EE_P(F)}
 \prod_{v \in P \setminus D}
  \frac{ 1 }
       { \vecta \langle H_P(v) \rangle },
\end{multline*}
where the summation on the left hand side is taken over all increasing sequences 
$F = F_0 \subsetneq F_1 \subsetneq \cdots \subsetneq F_n = P$ of order filters of length $n$, 
and $\vecta \langle D \rangle = \sum_{p \in D} a_{c(p)}$ for a subset $D \subset P$.
Such increasing sequences of order filters are in one-to-one correspondence 
with linear extensions $(q_1, \cdots, q_n)$ of $P \setminus F$ by the relation
$$
F_k = F \cup \{ q_n, \cdots, q_{n-k+1} \}
\quad(0 \le k \le n).
$$
Hence we obtain the desired result.
\end{proof}

\begin{rema}
Corollary~\ref{cor:main2} can be proved by using the Billey formula \cite[Theorem~4]{B} 
and the Chevalley formula \cite[Theorem~11.1.7 and Corollary~11.3.17]{Ku} 
for the equivariant cohomology along the same line as Theorem~\ref{thm:main1}.
\end{rema}

By specializing $z_i = q$ for all $i \in I$ in \eqref{eq:main1}, 
and $a_i = 1$ for all $i \in I$ in \eqref{eq:main2}, we obtain

\begin{coro}
\label{cor:main3}
Let $P$ be a $d$-complete poset and $F$ an order filter of $P$.
We define the hook length $h_P(u)$ at $u \in P$ as follows:
\begin{enumerate}
\item[(i)]
If $u$ is not the top of any $d_k$-interval, then we define $h_P(u) = \# \{ w \in P : w \le u \}$.
\item[(ii)]
If $u$ is the top of a $d_k$-interval $[v,u]$, then we define 
$h_P(u) = h_P(x) + h_P(y) - h_P(v)$, 
where $x$ and $y$ are the sides of $[v,u]$.
\end{enumerate}
Then we have
\begin{enumerate}
\item[(a)]
The univariate generating function of $(P \setminus F)$-partitions is given by
\begin{equation}
\label{eq:main31}
\sum_{\sigma \in \AP(P \setminus F)}
 q^{|\sigma|}
 =
\sum_{D \in \EE_P(F)}
 \frac{ \prod_{v \in B(D)} q^{h_P(v)} }
      { \prod_{v \in P \setminus D} ( 1 - q^{h_P(v)} ) }.
\end{equation}
\item[(b)]
The number of linear extensions of $P \setminus F$ is given by
\begin{equation}
\label{eq:main32}
n! \sum_{D \in \EE_P(F)} \prod_{v \in P \setminus D} \frac{1}{h_P(v)},
\end{equation}
where $n = \# (P \setminus F)$.
\end{enumerate}
\end{coro}

If $P = D(\lambda)$ and $F = D(\mu)$ are shapes corresponding to partitions $\lambda \supset \mu$,
Equations \eqref{eq:main31} and \eqref{eq:main32} reduce to 
Morales--Pak--Panova's $q$-hook formula \cite[Corollary~6.17]{MPP} and Naruse's hook formula \cite{Naruse} 
respectively.
The trace generating function of revers plane partitions of skew shape \cite[Corollary~6.20]{MPP}
is obtained from Theorem~\ref{thm:main1} by specializing
$$
z_i = \begin{cases}
 tq &\text{if $i$ is the color of the maximum element of $D(\lambda)$,} \\
 q &\text{otherwise.}
\end{cases}
$$

\begin{rema}
\label{rem:heap3}
Theorem~\ref{thm:main1} and its corollaries hold for heaps $H(w)$ associated 
to dominant minuscule elements $w$ in any symmetrizable Kac--Moody Weyl groups, 
after suitable modifications are made.
See Remarks~\ref{rem:heap1} and \ref{rem:heap2}.
\end{rema}

\subsection{%
Example
}

\begin{exam}
\label{ex:dc}
Let $P = S(3,2,1) \supset F = S(2)$ be the shifted shapes corresponding 
to strict partitions $(3,2,1)$ and $(2)$.
If we regard $P$ as a $d$-complete poset with a $d$-complete coloring $c : P \to \{ 0, 0', 1, 2 \}$ 
given in Example~\ref{ex:shifted_hc}, 
then the hook monomials in 
${{\vectz}} = (z_0, z_{0'}, z_1, z_2)$ are given by
\begin{alignat*}{3}
{{\vectz}}[ H_P(1,1) ] &= z_0 z_{0'} z_1^2 z_2,
&\quad
{{\vectz}}[ H_P(1,2) ] &= z_0 z_{0'} z_1 z_2,
&\quad
{{\vectz}}[ H_P(1,3) ] &= z_0 z_1 z_2,
\\
&
&
{{\vectz}}[ H_P(2,2) ] &= z_0 z_{0'} z_1,
&\quad
{{\vectz}}[ H_P(2,3) ] &= z_0 z_1,
\\
&
&
&
&{{\vectz}}[ H_P(3,3) ] &= z_0.
\end{alignat*}
Since we have
$$
\setlength{\unitlength}{1.1pt}
\mathcal{E}_P(F)
 =
\left\{
\raisebox{-14pt}{
\begin{picture}(30,30)
\fboxsep=0mm
\put(0,20){\colorbox[gray]{0.7}{\makebox(10,10){}}}
\put(10,20){\colorbox[gray]{0.7}{\makebox(10,10){}}}
\put(0,30){\line(1,0){30}}
\put(0,20){\line(1,0){30}}
\put(10,10){\line(1,0){20}}
\put(20,0){\line(1,0){10}}
\put(0,20){\line(0,1){10}}
\put(10,10){\line(0,1){20}}
\put(20,0){\line(0,1){30}}
\put(30,0){\line(0,1){30}}
\end{picture}
},
\ 
\raisebox{-14pt}{
\begin{picture}(30,30)
\fboxsep=0mm
\put(0,20){\colorbox[gray]{0.7}{\makebox(10,10){}}}
\put(20,10){\colorbox[gray]{0.7}{\makebox(10,10){}}}
\put(0,30){\line(1,0){30}}
\put(0,20){\line(1,0){30}}
\put(10,10){\line(1,0){20}}
\put(20,0){\line(1,0){10}}
\put(0,20){\line(0,1){10}}
\put(10,10){\line(0,1){20}}
\put(20,0){\line(0,1){30}}
\put(30,0){\line(0,1){30}}
\put(10,20){\makebox(10,10){$\times$}}
\end{picture}
}
\right\},
$$
we apply Theorem~\ref{thm:main1} to obtain
\begin{align}
\label{eq:ex1}
\sum_{\pi \in \AP(P \setminus F)} {\vectz}^\pi
 &=
\frac{ 1 }
     { (1 - z_0 z_1 z_2) (1 - z_0 z_{0'} z_1) (1 - z_0 z_1) (1 - z_0) }
\\
 &\quad
+
\frac{ z_0 z_{0'} z_1 z_2 }
     { (1 - z_0 z_{0'} z_1 z_2) (1 - z_0 z_1 z_2) (1 - z_0 z_{0'} z_1) (1 - z_0) }
\notag
\\
 &=
\frac{ 1 - z_0^2 z_{0'} z_1^2 z_2 }
     { (1 - z_0 z_{0'} z_1 z_2) (1 - z_0 z_1 z_2) (1 - z_0 z_{0'} z_1) (1 - z_0 z_1) (1 - z_0) },
\notag
\end{align}
where
$$
{\vectz}^\pi
 =
z_0^{\pi(1,1) + \pi(3,3)} z_{0'}^{\pi(2,2)} z_1^{\pi(1,2) + \pi(2,3)} z_2^{\pi(1,3)}.
$$
\end{exam}

\begin{exam}
\label{ex:heap}
Let $P = S(3,2,1) \supset F = S(2)$ be the same as in Example~\ref{ex:dc}.
If we regard $P$ as the heap $H(w_{(3,2,1)})$ for the Weyl group of type $B_3$ 
(see Example~\ref{ex:shifted-B}),
then the hook monomials in $\overline{\vectz} = (z_0, z_1, z_2)$ are given by
\begin{alignat*}{3}
\overline{\vectz}[ H'_P(1,1) ] &= z_0 z_1 z_2,
&\quad
\overline{\vectz}[ H'_P(1,2) ] &= z_0^2 z_1^2 z_2,
&\quad
\overline{\vectz}[ H'_P(1,3) ] &= z_0^2 z_1 z_2,
\\
&
&
\overline{\vectz}[ H'_P(2,2) ] &= z_0 z_1,
&\quad
\overline{\vectz}[ H'_P(2,3) ] &= z_0^2 z_1,
\\
&
&
&
&
\overline{\vectz}[ H'_P(3,3) ] &= z_0.
\end{alignat*}
Since we have
$$
\setlength{\unitlength}{1.1pt}
\mathcal{E}_P(F)
 =
\left\{
\raisebox{-14pt}{
\begin{picture}(30,30)
\fboxsep=0mm
\put(0,20){\colorbox[gray]{0.7}{\makebox(10,10){}}}
\put(10,20){\colorbox[gray]{0.7}{\makebox(10,10){}}}
\put(0,30){\line(1,0){30}}
\put(0,20){\line(1,0){30}}
\put(10,10){\line(1,0){20}}
\put(20,0){\line(1,0){10}}
\put(0,20){\line(0,1){10}}
\put(10,10){\line(0,1){20}}
\put(20,0){\line(0,1){30}}
\put(30,0){\line(0,1){30}}
\end{picture}
},
\ 
\raisebox{-14pt}{
\begin{picture}(30,30)
\fboxsep=0mm
\put(0,20){\colorbox[gray]{0.7}{\makebox(10,10){}}}
\put(20,10){\colorbox[gray]{0.7}{\makebox(10,10){}}}
\put(0,30){\line(1,0){30}}
\put(0,20){\line(1,0){30}}
\put(10,10){\line(1,0){20}}
\put(20,0){\line(1,0){10}}
\put(0,20){\line(0,1){10}}
\put(10,10){\line(0,1){20}}
\put(20,0){\line(0,1){30}}
\put(30,0){\line(0,1){30}}
\put(10,20){\makebox(10,10){$\times$}}
\end{picture}
},
\ 
\raisebox{-14pt}{
\begin{picture}(30,30)
\fboxsep=0mm
\put(10,10){\colorbox[gray]{0.7}{\makebox(10,10){}}}
\put(20,10){\colorbox[gray]{0.7}{\makebox(10,10){}}}
\put(0,30){\line(1,0){30}}
\put(0,20){\line(1,0){30}}
\put(10,10){\line(1,0){20}}
\put(20,0){\line(1,0){10}}
\put(0,20){\line(0,1){10}}
\put(10,10){\line(0,1){20}}
\put(20,0){\line(0,1){30}}
\put(30,0){\line(0,1){30}}
\put(0,20){\makebox(10,10){$\times$}}
\end{picture}
}
\right\},
$$
we apply a heap version of Theorem~\ref{thm:main1} to obtain
\begin{align}
\label{eq:ex2}
\sum_{\pi \in \AP(P \setminus F)} \overline{\vectz}^\pi
 &=
\frac{ 1 }
     { (1 - z_0^2 z_1 z_2) (1 - z_0 z_1) (1 - z_0^2 z_1) (1 - z_0) }
\\
 &\quad
+
\frac{ z_0^2 z_1^2 z_2 }
     { (1 - z_0^2 z_1^2 z_2) (1 - z_0^2 z_1 z_2) (1 - z_0 z_1) (1 - z_0) }
\notag
\\
 &\quad
+
\frac{ z_0 z_1 z_2 }
     { (1 - z_0 z_1 z_2) (1 - z_0^2 z_1^2 z_2) (1 - z_0^2 z_1 z_2) (1 - z_0) }
\notag
\\
 &=
\frac{ 1 - z_0^3 z_1^2 z_2 }
     { (1 - z_0 z_1 z_2) (1 - z_0^2 z_1 z_2) (1 - z_0 z_1) (1 - z_0^2 z_1) (1 - z_0) },
\notag
\end{align}
where
$$
\overline{\vectz}^\pi
 =
z_0^{\pi(1,1) + \pi(2,2) + \pi(3,3)} z_1^{\pi(1,2) + \pi(2,3)} z_2^{\pi(1,3)}.
$$
Note that Equation \eqref{eq:ex2} is obtained from \eqref{eq:ex1} by putting $z_{0'} = z_0$.
\end{exam}

\longthanks{
The authors are grateful to the two anonymous referees for their careful reading of our manuscript 
and their valuable comments and suggestions, which helped us to improve the exposition of this paper.
The first author was partially supported by the JSPS Grants-in-Aid for Scientific Research No.~16H03921.
The second author was partially supported by the JSPS Grants-in-Aid for Scientific Research No.~15K13425, 
and gratefully acknowledges the support and hospitality 
of the Erwin Schr\"odinger International Institute for Mathematics and Physics, 
where part of this work was carried out.
}

\nocite{*}
\bibliographystyle{amsplain-ac}
\bibliography{skewhook_alco}

\end{document}